\documentclass[final,1p,times]{elsarticle}

\usepackage{enumitem}
  \setlist{nosep} %

\usepackage[makeroom]{cancel}
\usepackage{amssymb,amsthm,mathtools,xfrac}
\usepackage{stmaryrd}
\usepackage{verbatim}
\usepackage{soul}
\usepackage{comment}
\usepackage{color}
\usepackage{cleveref}
\usepackage{subcaption}

\usepackage{algorithm,algorithmicx}
\usepackage[noend]{algpseudocode}

\newcommand*{\dif}{\mathop{}\!\mathrm{d}}

\usepackage{xcolor}

\newcommand{\tcb}[1]{\textcolor{black}{#1}}

\journal{Journal of Computational Physics}

\usepackage[mathlines]{lineno}
\newcommand*\patchAmsMathEnvironmentForLineno[1]{%
  \expandafter\let\csname old#1\expandafter\endcsname\csname #1\endcsname
  \expandafter\let\csname oldend#1\expandafter\endcsname\csname end#1\endcsname
  \renewenvironment{#1}%
     {\linenomath\csname old#1\endcsname}%
     {\csname oldend#1\endcsname\endlinenomath}}%
\newcommand*\patchBothAmsMathEnvironmentsForLineno[1]{%
  \patchAmsMathEnvironmentForLineno{#1}%
  \patchAmsMathEnvironmentForLineno{#1*}}%
\AtBeginDocument{%
\patchBothAmsMathEnvironmentsForLineno{equation}%
\patchBothAmsMathEnvironmentsForLineno{align}%
\patchBothAmsMathEnvironmentsForLineno{flalign}%
\patchBothAmsMathEnvironmentsForLineno{alignat}%
\patchBothAmsMathEnvironmentsForLineno{gather}%
\patchBothAmsMathEnvironmentsForLineno{multline}%
}

\begin{document}

\begin{frontmatter}

\title{Implicit-explicit Runge-Kutta for radiation hydrodynamics I: gray diffusion}
\author[LANL]{Ben S. Southworth}
\cortext[mycorrespondingauthor]{Corresponding author}
\ead{southworth@lanl.gov}
\author[LANL]{HyeongKae Park}
\author[LANL]{Svetlana Tokareva}
\author[Microsoft]{Marc Charest}

\address[LANL]{Theoretical Division, Los Alamos National Laboratory, P.O. Box 1663, Los Alamos, NM 87545 U.S.}

\begin{abstract}
Radiation hydrodynamics are a challenging multiscale and multiphysics set of equations. To capture the relevant physics of interest, one typically must time step on the hydrodynamics timescale, making explicit integration the obvious choice. On the other hand, the coupled radiation equations have a scaling such that implicit integration is effectively necessary in non-relativistic regimes. A first-order Lie-Trotter-like operator split is the most common time integration scheme used in practice, alternating between an explicit hydrodynamics step and an implicit radiation solve and energy deposition step. However, such a scheme is limited to first-order accuracy, and nonlinear coupling between the radiation and hydrodynamics equations makes a more general additive partitioning of the equations non-trivial. Here, we develop a new formulation and partitioning of radiation hydrodynamics with gray diffusion that allows us to apply (linearly) implicit-explicit Runge-Kutta time integration schemes. We prove conservation of total energy in the new framework, and demonstrate 2nd-order convergence in time on multiple radiative shock problems, achieving error 3--5 orders of magnitude smaller than the first-order Lie-Trotter operator split \emph{at the hydrodynamic CFL}, even when Lie-Trotter applies a 3rd-order TVD Runge-Kutta scheme to the hydrodynamics equations.
\end{abstract}

\end{frontmatter}

\section{Background}\label{sec:intro}

We are interested in solving the non-relativistic radiation-hydrodynamics equations with a grey radiation diffusion equation of the form:
\begin{subequations}\label{eq:rad-hydro}
\begin{align}
    \frac{\partial \rho}{\partial t} &=-  \nabla \cdot \left(\rho \mathbf{u}\right), \label{eq:rad-hydro-rho} \\
    \frac{\partial \rho \mathbf{u}}{\partial t} &=- \nabla \cdot \left[\rho \mathbf{u}\otimes\mathbf{u} + (p + {p_r})\mathbf{I}\right], \label{eq:rad-hydro-rhou}  \\
    \frac{\partial \rho e_t}{\partial t} &= -\nabla \cdot \left[ \left(\rho e_t + (p+{p_r})\right) \mathbf{u}\right]  + {p_r\nabla\cdot \mathbf{u}} + {\sigma_E cE_r} -{\sigma_p acT^4},\label{eq:rad-hydro-et}\\
    \frac{\partial E_r}{\partial t} &= -\nabla \cdot \left(E_r \mathbf{u}\right) - {p_r \nabla \cdot \mathbf{u}} + \nabla \cdot D \nabla E_r - {\sigma_E c E_r}+ {\sigma_p acT^4}, \label{eq:rad-hydro-Er}
\end{align}
\end{subequations}
where, $c$ is the speed of light, $\rho,~\mathbf{u},~e_t$ are the material density, velocity, and specific (hydro) total energy, respectively, $p$ and \tcb{$p_r=\frac{E_r}{3}$} are the hydro and radiation pressure, and $E_r$ is the radiation energy density. Total energy is defined as $E = \rho e_t + E_r$, and summing Eqs. \eqref{eq:rad-hydro-et} and \eqref{eq:rad-hydro-Er} yields the total energy equation,\footnote{Note, it is important to express \eqref{eq:rad-hydro-et} using $p_r\nabla\cdot \mathbf{u}$ (rather than simplifying with the chain rule applied to $\nabla\cdot p_r\mathbf{u}$) so that summing the \emph{discrete} form of \eqref{eq:rad-hydro-et} and \eqref{eq:rad-hydro-Er} maintains the cancellation of this term.}
\begin{equation}
        \frac{\partial E}{\partial t} = - \nabla \cdot \left[(E+p+p_r)\mathbf{u}\right] + \nabla \cdot D\nabla E_r. \label{eq:rad-hydro-E} 
\end{equation}
Specific internal energy is given by $e_i= e_t - u^2/2$, where (taking $\mathbf{u}$ to be a column vector) $u^2 = \|\mathbf{u}\|^2 = \mathbf{u}^T\mathbf{u}$. By a complete equation of state (EOS), opacities ($\sigma_E$ and $\sigma_p$), temperature, pressure, and the diffusion coefficient $D$ are nonlinear functions of energy and density.
 
In many contexts, to accurately capture the physics of interest, we must time step on the order of the hydrodynamic timescale, which also roughly corresponds with the hydrodynamics stability constraint for explicit time stepping. Thus, because we are already stepping within the explicit stability limit, the hydrodynamics equations will be treated explicitly to avoid the cost of implicit solves. However, the radiation energy and coupling to temperature propagate on the order of the speed of light. The resulting extreme limitation on timestep size for explicit integration leads to radiation almost universally being treated implicitly. The most common approach for the coupled system in \eqref{eq:rad-hydro} is then a variation on a first-order Lie-Trotter operator-splitting, combining an explicit hydrodynamics step with implicit radiation solves. To obtain second-order accuracy in time (which is not feasible with a Lie-Trotter operator split), some work has iterated between explicit hydrodynamics evaluations and implicit radiation solves in a single time step \cite{Kadioglu.2010,Bolding.2017} to resolve the nonlinear physics coupling in EOS and opacities to higher than first (linearized) order. A second-order Strang operator splitting \cite{strang1968construction} has also been applied to radiation hydrodynamics, e.g., \cite{Haines.2022,fuksman2021two}, but Strang splittings are known to suffer from order reduction, particularly for stiff problems (see, e.g., \cite{einkemmer2018efficient,sportisse2000analysis} in general and \cite{zingale2019improved} in the context of hydrodynamics). Recent work has also developed higher-order integration schemes from spectral deferred corrections (SDC) for reactive hydrodynamics \cite{zingale2019improved,zingale2022improved}, which are problem specific, but may provide a different tractable direction to achieve higher order for more complex hydrodynamics problems.

This paper introduces a general framework for second-order (and higher) time integration schemes for (non-relativistic) radiation hydrodynamics with gray diffusion, without requiring nonlinear iterations between the hydrodynamics and radiation equations, while also providing more rigorous and robust convergence properties than obtained through operator split. Broadly, this is achieved by moving away from a traditional operator split and instead integrating in time using carefully chosen Runge-Kutta methods. Runge Kutta methods and their additive and partitioned variations \cite{Rice.1960,hofer1976partially,hairer1981order, Kennedy.2003tv4} are widely used in numerical PDEs, and offer many practical advantages, including high-order accuracy, built in error-estimators, and good stability properties \cite{Hairer.1996ci}. However, the application of additive type Runge-Kutta methods to radiation hydrodynamics and similar multiphysics problems is not straightforward, because the physics do not always have a clear separation between stiff and non-stiff components. For example, material properties in radiation hydrodynamics are often linearized, and thus operators such as absorption, $\sigma_E cE_r$, cannot be naturally categorized as just implicit or explicit. In particular, the speed-of-light scaling makes the term too stiff to treat explicitly in $E_r$, but implicit treatment of $\sigma_E$ introduces nonlinear coupling to the temperature and hydrodynamics equations, significantly complicating an implicit solve. IMEX-RK methods have recently been applied to relativistic radiation hydrodynamics \cite{Fuksman.2019}, treating the stiff (and local) source terms implicitly; however, in the relativistic setting the explicit stability limit of hydrodynamics and radiation are comparable, and such an approach is not applicable in the non-relativistic case. We address these challenges by following the linearly implicit-explicit (LIMEX) approach of \cite{Boscarino.2015,Boscarino.2016,Boscarino.20168e8}, coupled with additive Runge Kutta (ARK) methods and a carefully derived ``partition'' of \eqref{eq:rad-hydro}. Such methods have recently demonstrated robust convergence applied to variations in the compressible Euler equations \cite{boscarino2018all,avgerinos2019linearly}, and prove key to addressing the coupling between radiation and hydrodynamics.

The paper proceeds as follows. \Cref{sec:time} details the LIMEX formulation for partitioned integration, providing background on general RK and ARK methods, the specific form of LIMEX-RK from \cite{Boscarino.2016}, and an additive LIMEX-RK (ALIMEX-RK) class of methods that proves useful to frame a Lie-Trotter operator split in the context of ARK methods. The partitioned integration of \eqref{eq:rad-hydro} is then discussed in \Cref{sec:hydro}. We begin by expressing the commonly used operator split integration as a single unified time step applied to \eqref{eq:rad-hydro}, and use this to derive ``partitions'' of \eqref{eq:rad-hydro} suitable for general LIMEX and ALIMEX methods in \Cref{sec:hydro:rk}. A key part of this is developing a consistent closure for the temperature that avoids nonlinearly coupling the EOS equations into the implicit radiation solve. In \Cref{sec:hydro:rk:alimex}, a first-order ALIMEX-RK method for integrating \eqref{eq:rad-hydro} is developed that is almost equivalent to the standard operator split (as it turns out, the standard operator split cannot be exactly posed as a standard ODE integrator due to inconsistent evaluation of temperature in the hydrodynamics vs. radiation equations). Conservation of energy under the proposed schemes is proven in \Cref{sec:energy}, and numerical results are provided in \Cref{sec:results} that demonstrate 2nd-order convergence in time on 1d radiative shock problems using the new methods. In addition to achieving 2nd-order accuracy in time, a major benefit of the new methods is a very small error constant at the hydrodynamics timescale; for a mach-3 radiative shock problem, the proposed integrators achieve error in hydrodynamics and radiation variables 3--4 orders of magnitude smaller than that achieved by a standard operator split, even when the operator split uses an inner 3rd order explicit RK method.

\section{Linearly-implicit time integration}\label{sec:time}

Suppose we have discretized our system of PDEs \eqref{eq:rad-hydro} in space, resulting in a nonlinear set of autonomous ODEs in time:
\begin{equation}\label{eq:ode}
    \frac{\partial y}{\partial t} = N(y)
\end{equation}
A widely used approach in physics simulations is to define an additive splitting of the operator into stiff and non-stiff components, $N(y) = N_I(y) + N_E(y)$, respectively, and proceed to treat $N_I(y)$ implicitly and $N_E(y)$ explicitly. As discussed in \Cref{sec:intro}, in radiation hydrodynamics the physics do not have a clear separation between stiff and non-stiff components. To mitigate this issue, we follow the linearly implicit-explicit (LIMEX) approach of \cite{Boscarino.2015,Boscarino.2016,Boscarino.20168e8}, instead introducing an auxiliary variable $y^*$ and defining the component partitioned set of ODEs
\begin{equation}\label{eq:ode-part}
    \frac{\partial y^*}{\partial t} = N(y^*,y),\hspace{3ex}
    \frac{\partial y}{\partial t} = N(y^*,y),
\end{equation}
with $y(t_0) = y^*(t_0) = y_0$. {\color{black} Note that because the right-hand sides of $y$ and $y^*$ are identical, given by $N(y^*,y)$, integrating the coupled system exactly in time will yield identical solutions, $y(t) = y^*(t)$, which are also identical to exact integration of the original equation \eqref{eq:ode}. By duplicating the equations, however, we are able to apply discrete partitioned integrators to the coupled system \eqref{eq:ode-part}, which are explicit in $y^*$ and implicit in $y$, regardless of nonlinear interaction between the two variables.} In doing this, we provide a flexibility and implementation very similar to standard Lie-Trotter operator splitting \emph{inside of the Runge Kutta stage solutions}, thereby incurring the benefits of Runge-Kutta integration, and avoiding limitations of standard operator splitting (most notably the limit to first-order accuracy for Lie-Trotter and a hard limit of 2nd-order accuracy for \emph{any} operator splitting).\footnote{It is worth pointing out that the recent paper \cite{Gonzalez-Pinto.2022} posed a Strang operator splitting in the context of generalized additive Runge Kutta (GARK) methods \cite{Sandu.2015}, but the standard formulation of GARK methods (or the framework in \cite{Gonzalez-Pinto.2022}) do not facilitate the non-additive nature of partitioning \eqref{eq:rad-hydro}.}

\subsection{Runge-Kutta methods}\label{sec:time:rk}

Standard Runge-Kutta methods are defined by the Butcher tableau
\[
\renewcommand\arraystretch{1.2}
\begin{array}
{c|c}
\mathbf{c} & A\\
\hline
& \mathbf{b}^T
\end{array}
\]
for coefficient matrix $A$, weights $\mathbf{b}^T$, and (local) quadrature nodes $\mathbf{c}$ (which do not come into play for autonomous ODEs). An $s$-stage Runge-Kutta method in stage-value formulation then takes the form
\begin{subequations}\label{eq:rk}
\begin{align}
    Y_{(i)} & = y^n + \Delta t \sum_{j=1}^s a_{ij} N(Y_{(j)}),\\
    y^{n+1} & = y^n + \Delta t \sum_{j=1}^s b_j N(Y_{(j)}).
\end{align}
\end{subequations}

Basic Runge-Kutta methods are often not viable for multiphysics systems due to the multiple timescales and stiffnesses observed; some physics will require implicit time integration due to minuscule explicit stability constraints (in the case of radiation, $\Delta t\lessapprox h/c$ for mesh spacing $h$ and speed of light $c$), while a fully coupled implicit solve over all physics is often intractably expensive. This led to the development of partitioned Runge Kutta (PRK) methods \cite{Rice.1960,hofer1976partially,hairer1981order} for ODEs that can be partitioned component wise, e.g., $y = [w;v]$ with $w'(t) = N_E(w(t),v(t)), v'(t) = N_I(w(t),v(t))$, and additive Runge Kutta (ARK) methods \cite{cooper1980additive,Cooper.1983sej} for ODEs that can be partitioned in an additive sense, e.g., $y' = N_E(y) + N_I(y)$. It turns out component-wise and additive partitions are equivalent classes, but one form typically provides a more natural representation for a given problem, and additive partitionings are more commonly considered.

Additive implicit-explicit (IMEX) Runge Kutta schemes were reintroduced and popularized in \cite{Ascher.1997}, where a number of schemes of various order and stability are presented and applied to an additively split advection-diffusion equation (see also \cite{Kennedy.2003tv4}). Suppose we split our operator into an explicit (non-stiff) and implicit (stiff) component,
\begin{align}\label{eq:imex}
    \frac{\partial y}{\partial t} = N(y) = N_E(y) + N_I(y).
\end{align}
We then define a combined (additive) IMEX-RK Butcher tableaux
\begin{equation}\label{eq:imex-tableaux}
\renewcommand\arraystretch{1.2}
\begin{array}
{c|c}
\tilde{\mathbf{c}} & \tilde{A}\\
\hline
&\tilde{\mathbf{b}}
\end{array}
,\hspace{3ex}
\begin{array}
{c|c}
\mathbf{c} & A\\
\hline
&\mathbf{b}
\end{array}
\end{equation}
where tilde-coefficients denote an explicit scheme ($\tilde{A}$ is strictly lower triangular, $\tilde{a}_{ij} = 0,~j \ge i$) and non-tilde coefficients represent a diagonally implicit Runge-Kutta (DIRK) scheme ($A$ is lower triangular). The classical IMEX Runge-Kutta scheme is then defined similar to \eqref{eq:rk} but in terms of two operators, coefficient sets, and weights:
\begin{subequations}\label{eq:imex-rk}
\begin{align}
    Y_{(i)} = y_n + \Delta t \sum_{j=1}^i \left(\tilde{a}_{ij} N_E(Y_{(j)}) + a_{ij} N_I(Y_{(j)})\right),\\
    y_{n+1} = y_n + \Delta t \sum_{j=1}^s \left( \tilde{b}_j N_E(Y_{(j)}) + b_j N_I(Y_{(j)})\right).
\end{align}
\end{subequations}
For example, by padding the tableaux \cite{Ascher.1997}, the simplest first-order IMEX-Euler method takes the form
\begin{equation*}
    y_{n+1} =  y_n + \Delta t N_E(y_n) + \Delta tN_I(y_{n+1}).
\end{equation*}
The order conditions for more general schemes follow that of partitioned RK methods \cite{Boscarino.2016}.

\subsection{Linearly implicit explicit Runge Kutta methods}\label{sec:time:ark}

In this paper we will split \eqref{eq:rad-hydro} foremost in a linearized sense, as in \eqref{eq:ode-part}. We propose a careful component partitioning to separate terms that are treated explicitly, linearly implicitly, and fully implicitly. The key is that we then solve a partitioned ODE of the form in \eqref{eq:ode-part}, wherein upon convergence $y=y^*$. Following \cite{Boscarino.2016}, by choosing an additive IMEX tableaux pair \eqref{eq:imex-tableaux} with $\mathbf{b} = \tilde{\mathbf{b}}$, a linearly implicit explicit Runge Kutta (LIMEX-RK) method for partitionings as in \eqref{eq:ode-part} takes the form
\begin{subequations}\label{eq:limex2}
\begin{align}
    Y^*_{(i)} &= y_n + \Delta t \sum_{j=1}^{i-1} \tilde{a}_{ij} N(Y^*_{(j)},Y_{(j)}), \\
    Y_{(i)}   &= y_n + \Delta t \sum_{j=1}^i {a_{ij}} N(Y^*_{(j)},Y_{(j)}),  \\
    y_{n+1} &= y_n + \Delta t \sum_{j=1}^s {b_{j}} N(Y^*_{(j)},Y_{(j)}). \label{eq:limex2-b}
\end{align}
\end{subequations}
For example, choosing the simplest IMEX-Euler scheme
\begin{equation}\label{eq:imex-euler}
\renewcommand\arraystretch{1.2}
\begin{array}
{c|c}
0 & 0\\
\hline
& 1
\end{array}
,\hspace{3ex}
\begin{array}
{c|c}
1 & 1\\
\hline
& 1
\end{array}
\end{equation}
the corresponding first-order LIMEX method is given by
\begin{equation}\label{eq:limex-o1}
    y_{n+1} = y_n + \Delta t N(y_n, y_{n+1}). 
\end{equation}
\tcb{Note that the purpose of choosing $\mathbf{b} = \tilde{\mathbf{b}}$ is so that the updated solutions $y_{n+1}$ and $y^*_{n+1}$ are identical summations \eqref{eq:limex2-b}, and thus we do not need to evaluate $y^*_{n+1}$ separately, or keep track of or choose which solution ($y_{n+1}$ or $y^*_{n+1}$) we move forward with. Since we are focused on autonomous problems, the quadrature points in time, $\mathbf{c}$ and $\tilde{\mathbf{c}}$, do not directly enter our integration scheme. However, it is worth pointing out that the schemes we have found to be most effective (see \Cref{sec:results:imex}) have distinct nodes $\mathbf{c}\neq\tilde{\mathbf{c}}$, and thus non-autonomous problems would take some extra care to avoid a $2\times$ increase in storing and evaluating stage vectors. The most efficient approach is to introduce time as a scalar auxiliary variable, and the quadrature node for a given stage is then where you evaluate this auxiliary variable, e.g. see \cite[Appendix]{Boscarino.2016}; such an approach requires care in implementation, but presents no technical or conceptual difficulties.} 

This \eqref{eq:limex2} is the class of methods we will use in practice for radiation hydrodynamics, and the IMEX tableaux we use are detailed in \Cref{sec:results:imex}. Note, from an implementation perspective, it can be useful to think of the nonlinear partition as
\begin{equation}\label{eq:limex-split}
    N(y^*,y) = N_E(y^*) + N_I(y^*,y),
\end{equation}
where $N_E$ is a purely explicit component (in our case hydrodynamics and material motion correction) and $N_I$ the semi-implicit component (radiation and energy deposition).

\subsection{Additive linearly implicit explicit Runge Kutta methods}\label{sec:time:ark2}

For the purpose of putting the standard radiation hydrodynamics operator split in a more general Runge Kutta framework, we also introduce a three-way partitioned additive linearly implicit explicit Runge Kutta (ALIMEX-RK) method (see \emph{Acknowledgements}). We use the semi-implicit-explicit partition in \eqref{eq:limex-split}, copy the equations as in the LIMEX setting, but now apply three different RK tableaux to the resulting system:
\begin{align}\label{eq:rad-hydro-part3}
    \frac{d}{dt}y^* & = \underbrace{N_I(y^*,y)}_{\hat{A}} + \underbrace{N_E(y^*)}_{\tilde{A}}, \hspace{3ex}
    \frac{d}{dt}y = \underbrace{N_I(y^*,y)}_{A} + \underbrace{N_E(y^*)}_{\tilde{A}},
\end{align}
where the underbraces indicate which of the three Runge-Kutta tableaux from a three-way additive Runge-Kutta method are applied to each operator:
\begin{equation}\label{eq:tableau-ark3}
\renewcommand\arraystretch{1.2}
\begin{array}
{c|c}
\hat{\mathbf{c}} & \hat{A}\\
\hline
& \hat{\mathbf{b}}
\end{array}, \hspace{3ex}
\begin{array}
{c|c}
\tilde{\mathbf{c}} & \tilde{A}\\
\hline
& \tilde{\mathbf{b}}
\end{array}, \hspace{3ex}
\begin{array}
{c|c}
{\mathbf{c}} & {A}\\
\hline
& {\mathbf{b}}
\end{array}.
\end{equation}
In our case, we construct $\hat{A}$ and $\tilde{A}$ to be explicit (strictly lower triangular) and $A$ to be of DIRK type (lower triangular). This corresponds to treating $N_E(y^*)$ explicit in both equations, and $N_I(y^*,y)$ explicit in $y^*$ (with a potentially different explicit step than $N_E(\cdot)$) and implicit in $y$. Then, a time step takes the form
\begin{subequations}\label{eq:ark3}
\begin{align}
    Y_{(i)}^* & = y_n^* + \Delta t\sum_{j=1}^{i-1} \hat{a}_{ij}N_I(Y_{(j)}^*,Y_{(j)}) + \Delta t\sum_{j=1}^{i-1}\tilde{a}_{ij}N_E(Y_{(j)}^*), \\
    Y_{(i)} & = y_n + \Delta t\sum_{j=1}^{i} a_{ij}N_I(Y_{(j)}^*,Y_{(j)})
    + \Delta t\sum_{j=1}^{i} \tilde{a}_{ij} N_E(Y_{(j)}^*), \\
    y_{n+1}^* & = y_n^* + \Delta t\sum_{j=1}^{s} \hat{b}_jN_I(Y_{(j)}^*,Y_{(j)}) + \Delta t\sum_{j=1}^{s}\tilde{b}_j N_E(Y_{(j)}^*), \\
    y_{n+1} & = y_n + \Delta t\sum_{j=1}^{s} b_jN_I(Y_{(j)}^*,Y_{(j)}) + \Delta t\sum_{j=1}^{s} \tilde{b}_jN_E(Y_{(j)}^*).
\end{align}
\end{subequations}
In theory, if $\hat{\mathbf{b}} \neq \mathbf{b}$ this method requires storing and advancing both $y_{n+1}^*$ and $y_{n+1}$. However, in practice for a $p$th order method, both variables are order $\mathcal{O}(\Delta t^p)$, and we will only store and propagate $y$ to the next time step (because it includes the implicit solve in $N_I(y^*,y)$), that is, we let $y_{n+1}^* = y_{n+1}$.

Consider the simple first order ARK
\begin{equation}\label{eq:3ark}
\renewcommand\arraystretch{1.2}
\begin{array}
{c|c c}
0 & 0 & 0 \\
1 & 0 & 0 \\
\hline
& 1 & 0
\end{array}, \hspace{3ex}
\begin{array}
{c|c c}
0 & 0 & 0 \\
1 & 1 & 0 \\
\hline
& 1 & 0
\end{array}, \hspace{3ex}
\begin{array}
{c|c c}
0 & 0 & 0 \\
1 & 0 & 1 \\
\hline
& 0 & 1
\end{array}.
\end{equation}
Then a time step for \eqref{eq:ark3}, dropping $y_{n+1}^*$, takes the form
\begin{subequations}\label{eq:ark3-first-order}
\begin{align}
    Y_{(2)}^* & = y_n + \Delta tN_E(y_n), \\
    y_{n+1} & = Y_{(2)} = y_n + \Delta tN_I(Y_{(2)}^*,Y_{(2)})
    + \Delta t N_E(y_n).
\end{align}
\end{subequations}
Here, we first take an explicit step with $N_E(\cdot)$, then linearize selected components of the implicit equation $N_I(\cdot,\cdot)$ with the updated explicit solution, before solving the remaining implicit equation for the new solution. This is in contrast to the first order LIMEX approach in \eqref{eq:limex-o1}, where the implicit linearization occurs about the previous time step $y_n$. There, the explicit part of the operator (rather than being split off in an additive manner \eqref{eq:imex}) is presumed to be represented via the $*$ variable in $N(y^*,y)$, and the explicit step in time is taken simultaneously with the implicit solve.

We emphasize that we only introduce the three-partition ALIMEX-RK approach as a Runge--Kutta realization of the classical first-order operator split typically used in radiation hydrodynamics (see \Cref{sec:hydro:rk:alimex}). We do not see any clear benefits over a two-partition LIMEX method \eqref{eq:limex2}, and the three-partition introduces additional complications of order conditions coupling three tableaux.

\section{Partitioning for radiation hydrodynamics}\label{sec:hydro}

\subsection{Lie-Trotter operator splitting}\label{sec:hydro:split}

Traditionally (and still in most DOE multiphysics packages) radiation hydrodynamics is solved via operator splitting, because the absorption-emission term in the radiation diffusion equation is stiff, while the hydrodynamics equations are time-advanced in an explicit manner. The standard operator-split time-stepping for \eqref{eq:rad-hydro} can be described as follows:
\begin{enumerate}
    \item[1.] \textit{Hydrodynamics step:}
    \begin{align*}
        \frac{\rho^{n+1}-\rho^{n}}{\Delta t} & = -  \nabla \cdot \left(\rho^n     \mathbf{u}^n\right),\\
        \frac{\rho^{n+1} \mathbf{u}^{n+1}-\rho^{n} \mathbf{u}^{n}}{\Delta t} & = -\nabla \cdot \left[\rho^n \mathbf{u}^n\otimes\mathbf{u}^n + (p^n + p_r^n)\mathbf{I}\right], \\
        \frac{\rho^{n+1} e_t^{*}-\rho^n e_t^{n}}{\Delta t} & = - \nabla \cdot \left[(\rho^n e_t^n+p^n+p_r^n)\mathbf{u}^n\right] + p_r^n \nabla \cdot \mathbf{u}^n.
    \end{align*}
    Then, we evaluate an intermediate internal energy $e^*_i = e_t^*-\|\mathbf{u}^{n+1}\|^2/2$, and temperature via EOS, $T^* = f_T(\rho^{n+1},e_i^*)$.
    
    \item[1.5] \textit{Radiation material motion correction:}
    \begin{equation*}
        \frac{E_r^{*}-E_r^n}{\Delta t} = - \nabla \cdot \left(E_r^{n} \mathbf{u}^n\right) - p_r^n \nabla \cdot \mathbf{u}^{n}.
    \end{equation*}
    
    \item[2.] \textit{Radiation solve step:} Evaluate opacities and diffusion coefficients based on the most up-to-date hydrodynamics variables, and solve the following implicit equations for $E_r^{n+1}$ and $T^{n+1}$ (with $*$ variables fixed as constant): 
    \begin{align*}
        \frac{E_r^{n+1}-E_r^*}{\Delta t} - \nabla \cdot D^{*} \nabla E_r^{n+1}+  \sigma_E^* cE_r^{n+1}&=   \sigma_p^* ac(T^{n+1})^4 ,\\
        \rho^{n+1} c_v \frac{T^{n+1}-T^*}{\Delta t} + \sigma_p^* ac(T^{n+1})^{4} - \sigma_E^* cE_r^{n+1} &= 0.
    \end{align*}
    
    \item[2.5] \textit{Energy deposition step:}
    \begin{equation}\label{eq:split-dep}
        e_i^{n+1} = e_i^* + c_v\left(T^{n+1}-T^{*} \right).
    \end{equation}
    Because kinetic energy does not change in the radiation step, we have $e^{n+1}_t = e^{*}_t + c_v(T^{n+1}-T^*)$, which then yields
    \begin{equation*}
        \frac{\rho^{n+1}e_t^{n+1}-\rho^{n+1}e_t^*}{\Delta t} = \frac{\rho^{n+1} c_v(T^{n+1}-T^*)}{\Delta t}
        = \sigma_E^*c E^{n+1} - \sigma_P^* ac(T^{n+1})^4.
    \end{equation*}
\end{enumerate}
Summing the full set of discrete equations shown above, we get a temporal discretization of the form:
\begin{subequations}\label{eq:op-split}
\begin{align}
    \frac{\rho^{n+1}-\rho^{n}}{\Delta t} & = -\nabla \cdot \left(\rho^n \mathbf{u}^n\right) , \label{eq:op-split-rho}\\
    \frac{\rho^{n+1} \mathbf{u}^{n+1}-\rho^{n} \mathbf{u}^{n}}{\Delta t} & = - \nabla \cdot \left[\rho^n \mathbf{u}^n\otimes\mathbf{u}^n + (p^n + p_r^n)\mathbf{I}\right], \label{eq:op-split-rhou} \\
    \begin{split}
        \frac{\rho^{n+1} e_t^{n+1}-\rho^n e_t^{n}}{\Delta t} &= -\nabla \cdot \left[(\rho^n e_t^n+p^n+p_r^n)\mathbf{u}^n\right] + p_r^n \nabla \cdot \mathbf{u}^n + \\&\hspace{10ex}\sigma_E^* cE_r^{n+1} - \sigma_p^* ac(T^{n+1})^4,
    \end{split}\label{eq:op-split-rhoe}\\
    \begin{split}
    \frac{E_r^{n+1}-E_r^n}{\Delta t} & = -\nabla \cdot \left(E_r^{n} \mathbf{u}^n\right) - p_r^n \nabla \cdot \mathbf{u}^{n} + \nabla \cdot D^{*} \nabla E_r^{n+1} - \\&\hspace{10ex}\sigma_E^* c E_r^{n+1} + \sigma_p^* ac(T^{n+1})^4,
    \end{split}\label{eq:op-split-Er} \\
    \rho^{n+1} c_v \frac{T^{n+1}-T^*}{\Delta t} & = \sigma_E^* cE_r^{n+1} - \sigma_p^* ac(T^{n+1})^{4},\label{eq:op-split-T}
\end{align}
\end{subequations}
where $T^* = f_T(\rho^{n+1}, e_t^*-\|\mathbf{u}^{n+1}\|^2/2)$ for EOS $f_T(\rho,e_i)$ and hydro specific energy $e_t^*$ following the explicit hydrodynamics step.

\subsection{A closure in temperature}\label{sec:hydro:T}

A closure in pressure for \eqref{eq:rad-hydro} is defined via an equation of state (EOS), $p = f_p(\rho,e_i)$. One can also pose a closure in temperature via an EOS, $T = f_T(\rho, e_i)$; however, using such a closure in the radiation equation introduces nonlinear coupling to tabular EOS data, which is inefficient, difficult to linearize, and something that we want to avoid in practice. Instead, recall for differential $\dif$, in a closed thermodynamics system we have change in internal energy is proportional to change in temperature,
\begin{equation} \label{eq:balance0}
    \dif e_i = c_v \dif T 
        \hspace{3ex}\Longleftrightarrow\hspace{3ex}
    \frac{\partial e_i}{\partial T} = c_v
        \hspace{3ex}\Longleftrightarrow\hspace{3ex}
    \frac{\partial e_i}{\partial t } = c_v \frac{\partial T}{\partial t}.
\end{equation}

As a result, we formulate a time-dependent equation for temperature based on the requirement that change in internal energy is proportional to heat capacity times the change in temperature \eqref{eq:balance0}. Recalling that $e_i = e_t - \mathbf{u}^T\mathbf{u}/2$ and multiplying by $\rho$, we then have 
\begin{equation}\label{eq:temp-eq0}
  \rho c_v \tfrac{d}{dt} T = \rho\tfrac{d}{dt} e_t - \rho\tfrac{d}{dt}(\mathbf{u}^T\mathbf{u})/2.
\end{equation}
Note the identities $\tfrac{d}{dt}(\rho e_t) = \rho\tfrac{d}{dt} e_t+ e_t\tfrac{d}{dt}\rho$ and $\rho\tfrac{d}{dt}(\mathbf{u}^T\mathbf{u})/2 = \rho\mathbf{u}^T\tfrac{d}{dt}\mathbf{u} = \mathbf{u}^T\tfrac{d}{dt}(\rho\mathbf{u}) - \mathbf{u}^T\mathbf{u}\tfrac{d}{dt}\rho$.
Plugging into \eqref{eq:temp-eq0} yields 
\begin{equation}\label{eq:temp-eq1}
  \rho c_v \tfrac{d}{dt}T = \tfrac{d}{dt}(\rho e_t) -
    e_t \tfrac{d}{dt}\rho -
    \mathbf{u}^T\tfrac{d}{dt}( \rho\mathbf{u}) +
    \mathbf{u}^T\mathbf{u}\tfrac{d}{dt} \rho.
\end{equation}
Define 
\begin{subequations}
\begin{align}
    N_E^{\rho}(\rho,\mathbf{u}) & \coloneqq -\nabla \cdot \left(\rho     \mathbf{u}\right), \\
    N_E^{\rho\mathbf{u}}(\rho,\mathbf{u},e_t) & \coloneqq - \nabla \cdot \left[\rho \mathbf{u}\otimes\mathbf{u} + (p + p_r)\mathbf{I}\right], \\
    N_E^{\rho e_t}(\rho,\mathbf{u},e_t) & \coloneqq - \nabla \cdot \left[(\rho e_t+p+p_r)\mathbf{u}\right] + p_r \nabla \cdot \mathbf{u},
\end{align}
\end{subequations}
corresponding to the hydrodynamics portion of \eqref{eq:rad-hydro}, which will be treated explicitly. Noting that $\tfrac{d}{dt}(\rho e_t) = N_E^{\rho e_t} + \sigma_E cE_r - \sigma_p acT^4$, $\tfrac{d}{dt} (\rho\mathbf{u}) = N_E^{\rho\mathbf{u}}$, and $\tfrac{d}{dt} \rho = N_E^\rho$, we now modify the temperature equation \eqref{eq:temp-eq1} to take the form 
\begin{subequations}\label{eq:temp-final}
\begin{align}
    \tfrac{d}{dt} T & = \tfrac{1}{\rho c_v}\left( \sigma_E cE_r - \sigma_p acT^4 + 
        \mathcal{L}_T(\rho,\mathbf{u},e_t) \right), \hspace{3ex}\textnormal{where} \\
    \mathcal{L}_T(\rho,\mathbf{u},e_t) & = N_E^{\rho e_t}(\rho,\mathbf{u},e_t) - \mathbf{u}^TN_E^{\rho\mathbf{u}}(\rho,\mathbf{u},e_t) + (\mathbf{u}^T\mathbf{u}-e_t)N_E^\rho(\rho,\mathbf{u}).
\end{align}
\end{subequations}

What we have here is effectively an internal energy equation written in terms of temperature, for the purposes of coupling with the implicit radiation solve. For consistency purposes, we have constructed things in such a way that we can evaluate the nonlinear temperature equation directly using other solution variables and the corresponding spatial discretizations that arise for each time derivative. Note, \emph{at the beginning of each time step, we define the current temperature directly via evaluating an EOS at the given hydrodynamics variables.} Because we are utilizing single-step integration methods, this is not relevant for the time integration schemes, and simply ensure greater fidelity to the temperature EOS.

\subsection{Runge-Kutta partitions}\label{sec:hydro:rk}

When we partition equations for an additive Runge Kutta or multistep type method, we do not form two time-dependent subsystems and take successive time steps, rather we construct of an additive partition of the right-hand side of the system of ODEs. As a result, such a partition is subtly different than that used in operator splitting. Here, we express the full spatial operator as the sum of two block operators over our variables $[\rho, \rho \mathbf{u}, \rho e_t, E_r, T]$:
\begin{enumerate}
    \item \textit{Explicit hydrodynamics and radiation MMC partition:}
    \begin{equation}\label{eq:rk-part-exp}
        \begin{bmatrix}
            N_E^{\rho} \\ N_E^{\rho\mathbf{u}} \\ N_E^{\rho e_t} \\ N_E^{E_r} \\ N_E^{T}
        \end{bmatrix}
        \coloneqq
        \begin{bmatrix}
            -\nabla \cdot \left(\rho     \mathbf{u}\right) \\
            - \nabla \cdot \left[\rho \mathbf{u}\otimes\mathbf{u} + (p + p_r)\mathbf{I}\right] \\
            - \nabla \cdot \left[(\rho e_t+p+p_r)\mathbf{u}\right] + p_r \nabla \cdot \mathbf{u} \\
            - \nabla \cdot \left(E_r \mathbf{u}\right) - p_r \nabla \cdot \mathbf{u} \\ 
            0
        \end{bmatrix}
    \end{equation}
    
    \item \textit{Radiation solve and energy deposition step:} This is again the subtle difference between operator split and partitioned integration. In general partitioned integration, we need the partitions to sum to the original operator, i.e., we cannot use the strong substitution of $\tfrac{d}{dt} \rho = 0$ to map $\rho \tfrac{d}{dt} e_i\mapsto \tfrac{d}{dt} (\rho e_t)$ in one operator as used in the operator-split formulations. As a result, we derive a time-dependent temperature equation from the proportionality of change in internal energy to change in temperature as in \eqref{eq:temp-final}. From here, we arrive with the second part of our linearly implicit additive partition:
    \begin{equation}\label{eq:rk-part-imp}
        \begin{bmatrix}
            N_I^{\rho} \\ N_I^{\rho\mathbf{u}} \\ N_I^{\rho e_t} \\ N_I^{E_r} \\ N_I^{T}
        \end{bmatrix}
        =
        \begin{bmatrix}
            0 \\ 0 \\ \sigma_E^* cE_r - \sigma_p^* acT^4 \\
            \nabla \cdot D^* \nabla E_r -  \sigma_E^* cE_r +   \sigma_p^* acT^4 \\
            \left[ \sigma_E^* cE_r - \sigma_p^* acT^4 + \mathcal{L}_T(\rho^*,\mathbf{u}^*,e_t^*)\right] / (\rho^*c_v^*)
        \end{bmatrix}.
    \end{equation}
    Note, here $*$ indicates quantities that will be linearized via an auxiliary variable, as in \eqref{eq:ode-part}.
\end{enumerate}

\subsubsection{LIMEX-RK method}\label{sec:hydro:rk:limex}

Here we consider posing radiation hydrodynamics \eqref{eq:rad-hydro} in the LIMEX-RK form \eqref{eq:ode-part}. We define a linearly implicit explicit operator by summing the implicit \eqref{eq:rk-part-imp} and explicit \eqref{eq:rk-part-exp} operators expressed above, and evaluating all explicit quantities as $*$ variables. For ease of presentation, we separate out the fully explicit component and semi-implicit component as in \eqref{eq:limex-split}:
\begin{subequations}
\begin{align}
    N_E\left(
        \begin{bmatrix}
            \rho^* \\ \rho^* \mathbf{u}^* \\ \rho^* e_t^* \\ E_r^* \\ T^*
        \end{bmatrix}
    \right)
    & =
    \begin{bmatrix}
        -\nabla \cdot \left(\rho^*\mathbf{u}^*\right) \\
        - \nabla \cdot \left[\rho^* \mathbf{u}^*\otimes\mathbf{u}^* + (p^* + p_r^*)\mathbf{I}\right] \\
        - \nabla \cdot \left[(\rho^* e_t^*+p^*+p_r^*)\mathbf{u}^*\right] + p_r^* \nabla \cdot \mathbf{u}^* \\
        - \nabla \cdot \left(E_r^* \mathbf{u}^*\right) - p_r^* \nabla \cdot \mathbf{u}^* \\ 
        \mathcal{L}_T(\rho^*,\mathbf{u}^*,e_t^*) / (\rho^*c_v^*)
    \end{bmatrix}, \\
    N_I\left(
        \begin{bmatrix}
            \rho^* \\ \rho^* \mathbf{u}^* \\ \rho^* e_t^* \\ E_r^* \\ T^*
        \end{bmatrix},
        \begin{bmatrix}
            \rho \\ \rho \mathbf{u} \\ \rho e_t \\ E_r \\ T
        \end{bmatrix}
    \right)
    & =
    \begin{bmatrix}
        0 \\
        0 \\
        \sigma_E^* cE_r - \sigma_p^* acT^4 \\
        \nabla \cdot D^* \nabla E_r -  \sigma_E^* cE_r + \sigma_p^* acT^4\\ 
        \left[\sigma_E^* cE_r - \sigma_p^* acT^4\right] / (\rho^*c_v^*)
    \end{bmatrix}.
\end{align}
\end{subequations}
This provides a framework to apply general LIMEX-RK methods as derived in \cite{Boscarino.2016} and discussed in \Cref{sec:time:rk}.

The simplest LIMEX-RK method is the first order Euler method \eqref{eq:limex-o1}. For the radiation hydrodynamics splitting proposed here, this leads to the time integration scheme
\begin{subequations}\label{eq:ark-split1}
\begin{align}
    \frac{\rho^{n+1}-\rho^{n}}{\Delta t} & =-\nabla \cdot \left(\rho^n\mathbf{u}^n\right), \\
    \frac{\rho^{n+1} \mathbf{u}^{n+1}-\rho^{n} \mathbf{u}^{n}}{\Delta t} & = - \nabla \cdot \left[\rho^n\mathbf{u}^n\otimes\mathbf{u}^n + (p^n + p_r^n)\mathbf{I}\right], \\
    \begin{split}
    \frac{\rho^{n+1} e_t^{n+1}-\rho^n e_t^{n}}{\Delta t} & = - \nabla \cdot \left[(\rho^n e_t^n+p^n+p_r^n)\mathbf{u}\right] + p_r^n \nabla \cdot \mathbf{u}^n \\&\hspace{10ex}+ \sigma_E^n cE_r^{n+1} - \sigma_p^n ac(T^{n+1})^4,
    \end{split}\label{eq:ark-split1-hydroE}\\
    \begin{split}
    \frac{E_r^{n+1}-E_r^n}{\Delta t} & = - \nabla \cdot \left(E_r^n \mathbf{u}^n\right) - p_r^n \nabla \cdot \mathbf{u}^n + \nabla \cdot D^n \nabla E_r^{n+1} - \\&\hspace{10ex} \sigma_E^n cE_r^{n+1} +   \sigma_p^n ac(T^{n+1})^4,
    \end{split}\label{eq:ark-split1-radE}\\ 
   (\rho^{n+1}c_v) \frac{T^{n+1}-T^n}{\Delta t} & = \sigma_E^n cE_r^{n+1} - \sigma_p^n ac(T^{n+1})^4 + \mathcal{L}_T(\rho^n,\mathbf{u}^n,e_t^n).
\end{align}
\end{subequations}
This corresponds to an operator-split-like method, that linearizes implicit variables and hydrodynamic temperature update based upon the previous time step, rather than following the explicit hydrodynamics step as in \Cref{sec:hydro:split}. Higher order methods as introduced in \Cref{sec:results:imex} will be utilized in the numerical results. Recall from \Cref{sec:hydro:T}, at the beginning of each time step we define $T_n$ via an EOS evaluated at the current hydrodynamic variables. 

\subsubsection{ALIMEX-RK and equivalence to operator split}\label{sec:hydro:rk:alimex}

To reproduce the operator split that updates opacities and diffusion coefficients in the implicit radiation equation with the updated hydrodynamic variables, we need the additive LIMEX scheme from \eqref{eq:ark3-first-order}, coupled with an additive partition of the rad-hydro equations. Define a continuous time splitting, with $N_I$ dependent on $y$ and $y^*$ variables, as
\begin{align*}
    N_E(y) \coloneqq
    \begin{bmatrix}
        N_E^{\rho} \\ N_E^{\rho\mathbf{u}} \\ N_E^{\rho e_t} \\ N_E^{E_r} \\ N_E^{T}
    \end{bmatrix}
    & =
    \begin{bmatrix}
        -\nabla \cdot \left(\rho     \mathbf{u}\right) \\
        - \nabla \cdot \left[\rho \mathbf{u}\otimes\mathbf{u} + (p + p_r)\mathbf{I}\right] \\
        - \nabla \cdot \left[(\rho e_t+p+p_r)\mathbf{u}\right] + p_r \nabla \cdot \mathbf{u} \\
        - \nabla \cdot \left(E_r \mathbf{u}\right) - p_r \nabla \cdot \mathbf{u} \\ 
        0
    \end{bmatrix},\\
    N_I(y^*,y) \coloneqq
    \begin{bmatrix}
        N_I^{\rho} \\ N_I^{\rho\mathbf{u}} \\ N_I^{\rho e_t} \\ N_I^{E_r} \\ N_I^{T}
    \end{bmatrix}
    & =
    \begin{bmatrix}
        0 \\ 0 \\ \sigma_E^* cE_r - \sigma_p^* acT^4 \\
        \nabla \cdot D^* \nabla E_r -  \sigma_E^* cE_r +   \sigma_p acT^4 \\
        \left[\sigma_E^* cE_r - \sigma_p^* acT^4 + \mathcal{L}_T(\rho^*,\mathbf{u}^*,e_t^*)\right]/(\rho^*c_v^*)
    \end{bmatrix}.
\end{align*}
Plugging into \eqref{eq:ark3-first-order} leads to the updates
\begin{subequations}\label{eq:ark-split2}
\begin{align}
    \frac{\rho^{n+1}-\rho^{n}}{\Delta t} & =-\nabla \cdot \left(\rho^n\mathbf{u}^n\right), \\
    \frac{\rho^{n+1} \mathbf{u}^{n+1}-\rho^{n} \mathbf{u}^{n}}{\Delta t} & = - \nabla \cdot \left[\rho^n\mathbf{u}^n\otimes\mathbf{u}^n + (p^n + p_r^n)\mathbf{I}\right], \\
    \begin{split} 
    \frac{\rho^{n+1} e_t^{n+1}-\rho^n e_t^{n}}{\Delta t} & = - \nabla \cdot \left[(\rho^n e_t^n+p^n+p_r^n)\mathbf{u}\right] + p_r^n \nabla \cdot \mathbf{u}^n + \\&\hspace{10ex} \sigma_E^* cE_r^{n+1} - \sigma_p^* ac(T^{n+1})^4,
    \end{split} \\
    \begin{split} 
    \frac{E_r^{n+1}-E_r^n}{\Delta t} & = - \nabla \cdot \left(E_r^n \mathbf{u}^n\right) - p_r^n \nabla \cdot \mathbf{u}^n + \nabla \cdot D^* \nabla E_r^{n+1} - \\&\hspace{10ex}  \sigma_E^* cE_r^{n+1} +   \sigma_p^* acT^4,
    \end{split}\\ 
    (\rho^{n+1}c_v)\frac{T^{n+1}-T^n}{\Delta t} & = \sigma_E^* cE_r - \sigma_p^* ac(T^{n+1})^4 + \mathcal{L}_T(\rho^*,\mathbf{u}^*,e_t^*). \label{eq:ark-split2-T}
\end{align}
\end{subequations}
Here we almost exactly reproduce the Lie-Trotter operator splitting, with a slight modification in the temperature equation. In the operator split \eqref{eq:op-split}, we define the temperature equation to be solved implicitly with initial value $T^* = f_T(\rho^*,e_i^*)$. If we instead substitute
\begin{equation}\label{eq:approx-T*}
    T^* = T_n + \Delta t\mathcal{L}_T(\rho^*,\mathbf{u}^*,e_t^*) / (\rho^*c_v^*)
\end{equation}
into the operator split temperature equation \eqref{eq:op-split-T} or the ARK temperature equation \eqref{eq:ark-split2-T}, we arrive at identical methods in \eqref{eq:op-split} and \eqref{eq:ark-split2}.

For a more intuitive interpretation, recall that the time-dependent temperature equation is auxiliary in nature, in the sense that given other state variables we can always evaluate the temperature EOS. In operator split, we update the temperature value directly from the EOS following the hydrodynamics time step. This can be thought of as the first step in a LIMEX-type integration, where we exactly integrate changes in temperature due to changes in hydrodynamic variables. We then follow with a first order linearized implicit equation for the update in temperature arising from radiation energy (which in theory, negates any higher-order accuracy in temperature provided by evaluating the EOS directly following the hydrodynamics step). In the ARK setting, the substitution in \eqref{eq:approx-T*} instead corresponds to a \emph{first order} approximation of the temperature update following the explicit hydrodynamics step, in this case derived from the time-dependent temperature equations \eqref{eq:temp-final} rather than directly from the EOS. Analogous to the operator split, this is then followed with a first order implicit equation updating temperature with respect to radiation. Thus although the methods are not exactly equivalent, they are conceptually equivalent and the formal accuracy is first order in both cases.

Numerical results in \Cref{sec:results} demonstrate that for the three radiative shock problems considered, the Lie-Trotter operator split \eqref{eq:op-split} and first-order LIMEX-Euler \eqref{eq:ark-split1} provide near identical accuracy. The ALIMEX method in \eqref{eq:ark-split2} is not considered in practice, as we see no practical benefit of a three-way additive partition, and expect the solution to be almost identical to the operator split. 

\section{Conservation of energy}\label{sec:energy}

\subsection{Local conservation}

Here we provide a formal proof of local conservation of total energy $E = \rho e_t + E_r$ for the first-order time integration scheme \eqref{eq:ark-split1} in the newly proposed LIMEX-RK framework. Summing the equations \eqref{eq:ark-split1-hydroE} and \eqref{eq:ark-split1-radE} gives
\begin{equation}
\label{eq:enrg-total}
    \frac{E^{n+1} - E^n}{\Delta t} = \\ - \nabla \cdot \left[\left(\rho^n e_t^n+p^n+p_r^n + E_r^n\right) \mathbf{u}^n\right] + \nabla \cdot D^n \nabla E_r^{n+1}
\end{equation}
A fully discrete version of this equation must be consistent with the continuous equation for the total energy \eqref{eq:rad-hydro-E}. Assume that equations \eqref{eq:ark-split1-hydroE} and \eqref{eq:ark-split1-radE} have been discretized in space using a finite volume method as follows:
\begin{subequations}
\begin{align}
\label{eq:ark-split1-rho_et-discrete}
\begin{split}
 \frac{\overline{\rho e_t}^{n+1}-\overline{\rho e_t}^{n}}{\Delta t} |K| & = - \int\limits_{\partial K}\left[(\rho^n e_t^n+p^n+p_r^n)\mathbf{u}^n\right]\cdot\mathbf{n}\,dl + \\ &\hspace{5ex} \int\limits_{K} \left[p_r^n \nabla \cdot \mathbf{u}^n + \sigma_E^n cE_r^{n+1} - \sigma_p^n ac(T^{n+1})^4 \right]\,d\mathbf{x},
\end{split}\\
\label{eq:ark-split1-Er-discrete}
\begin{split}
 \frac{\overline{E_r}^{n+1}-\overline{E_r}^n}{\Delta t}|K| & = - \int\limits_{\partial K}\left(E_r^n \mathbf{u}^n\right)\cdot\mathbf{n}\,dl - \\ &\hspace{-15ex} \int\limits_{K} \left[p_r^n \nabla \cdot \mathbf{u}^n + \sigma_E^n cE_r^{n+1} - \sigma_p^n ac(T^{n+1})^4 \right]\,d\mathbf{x} + \int\limits_K \nabla \cdot D^n \nabla E_r^{n+1} \,d\mathbf{x},
 \end{split}
\end{align}
\end{subequations}
where $\overline{f}$ denotes the volume average on a given element $K$, and $|K|$ denotes the volume of an element. Note that for now we will not discuss the computation of the volume integrals and the diffusion term, it is  important, however, that the discretization of these terms is the same in \eqref{eq:ark-split1-rho_et-discrete} and \eqref{eq:ark-split1-Er-discrete}.

The surface integrals in \eqref{eq:ark-split1-rho_et-discrete} and \eqref{eq:ark-split1-Er-discrete} will be computed using numerical flux functions (due to discontinuities in the numerical solution at the cell interfaces). Denote the corresponding physical fluxes by $F_{\rho e_t} = (\rho e_t+p+p_r)\mathbf{u}$ and $F_{E_r} = E_r \mathbf{u}$, and the total energy flux by $F_E = (\rho e_t+p+p_r+E_r)\mathbf{u}$. There are many ways to compute the numerical flux function for a given physical flux $F(U)$, e.g. we could use local Lax-Friedrichs fluxes:
\begin{align}
    \hat{F}^{LxF} = \frac12(F(U_L) + F(U_R)) - \frac12\alpha(U_R-U_L),
\end{align}
where $U_L$ and $U_R$ are solutions in the neighboring mesh elements that share the interface $\partial K$.

Then the discrete scheme becomes
\begin{subequations}
\begin{align}
\label{eq:ark-split1-rho_et-discrete-1}
\begin{split}
 \frac{\overline{\rho e_t}^{n+1}-\overline{\rho e_t}^{n}}{\Delta t} |K| &= - \int\limits_{\partial K}\hat{F}_{\rho e_t}^n\cdot\mathbf{n}\,dl + \\& \int\limits_{K} \left[p_r^n \nabla \cdot \mathbf{u}^n + \sigma_E^n cE_r^{n+1} - \sigma_p^n ac(T^{n+1})^4 \right]\,d\mathbf{x},
\end{split} \\
\label{eq:ark-split1-Er-discrete-1}
\begin{split}
 \frac{\overline{E_r}^{n+1}-\overline{E_r}^n}{\Delta t}|K| & = - \int\limits_{\partial K}\hat{F}_{E_r}^n\cdot\mathbf{n}\,dl  - \\&\hspace{-15ex} \int\limits_{K} \left[p_r^n \nabla \cdot \mathbf{u}^n + \sigma_E^n cE_r^{n+1} - \sigma_p^n ac(T^{n+1})^4 \right]\,d\mathbf{x} + \int\limits_K \nabla \cdot D^n \nabla E_r^{n+1} \,d\mathbf{x},
\end{split}
\end{align}
\end{subequations}
and fully discrete version of \eqref{eq:enrg-total} is given by
\begin{equation}
\label{eq:enrg-total-discrete}
 \frac{\overline{E}^{n+1} - \overline{E}^n}{\Delta t} |K| = \\ - \int\limits_{\partial K} \hat{F}_{E}^n\cdot\mathbf{n}\,dl  + \int\limits_K \nabla \cdot D^n \nabla E_r^{n+1} \,d\mathbf{x}.
\end{equation}

Comparing equations \eqref{eq:ark-split1-rho_et-discrete}, \eqref{eq:ark-split1-Er-discrete}, \eqref{eq:enrg-total-discrete}, we can conclude that the scheme is locally conservative at the discrete level if the following consistency condition on the numerical fluxes is satisfied:
\begin{equation}
\label{eq:fluxes-conservation}
\hat{F}_{E}^n = \hat{F}_{\rho e_t}^n + \hat{F}_{E_r}^n,
\end{equation}
which can be guaranteed for the local Lax-Friedrichs flux as long as both $\hat{F}_{\rho e_t}^n$ and  $\hat{F}_{E_r}^n$ use the same viscosity coefficient $\alpha$. 

\subsection{Global conservation}

Local conservation as stated in \eqref{eq:fluxes-conservation} naturally yields global conservation of the total energy. Taking the summation over all elements $K$ in the discrete domain, we can write
\begin{multline}
\label{eq:enrg-global-conservation}
 \sum\limits_K|K|\overline{E}^{n+1} = \sum\limits_K|K|\overline{E}^n - {\Delta t}\left[ \sum\limits_K\int\limits_{\partial K} \hat{F}_{E}^n\cdot\mathbf{n}\,dl  - \sum\limits_K\int\limits_K \nabla \cdot D^n \nabla E_r^{n+1} \,d\mathbf{x}\right] = \\
 \sum\limits_K|K|\overline{E}^n - {\Delta t}\left[ \int\limits_{\Gamma_B} \hat{F}_{E}^n\cdot\mathbf{n}\,dl - \int\limits_{\Gamma_B} \widehat{(D^n \nabla E_r^{n+1})}\cdot\mathbf{n}\,dl \right],
\end{multline}
where $\Gamma_B$ is the physical boundary of the computational domain and $\widehat{(D^n \nabla E_r^{n+1})}$ is a numerical flux for the diffusion operator. Hence, the global conservation of energy is determined by the influx of energy and diffusion through the physical boundaries. In the absence of boundary flux, we have a perfect balance of discrete total energy, i.e. $\sum\limits_K|K|\overline{E}^{n+1} = \sum\limits_K|K|\overline{E}^{n}$.

\subsection{General LIMEX integration schemes}

Above, we established conditions for local and global conservation for the first-order LIMEX scheme in \eqref{eq:ark-split1}. To establish conservation for general LIMEX schemes, we simply point out that each stage of the general LIMEX scheme \eqref{eq:limex2} takes exactly the form of a scaled (by $b_j$) LIMEX-Euler step \eqref{eq:ark-split1}. Thus, again in the absence of boundary flux over each internal LIMEX-RK stage, we maintain conservation of energy when summing over all stages to get the new solution. Note that although conservation of linear invariants such as energy is trivial for explicit and implicit integration schemes \cite{Hairer.1996ci}, it becomes much more challenging for IMEX integration when the different partitions of the operator are evaluated at different stage vectors. We found the LIMEX approach to provide a more natural route to ensuring conservation for this problem than more common ARK-like methods. 

\section{Numerical results -- radiative shocks}\label{sec:results}

To test the newly proposed integration schemes, we consider multiple 1d radiative shock problems. The spatial discretization is finite-volume-based for the hydrodynamics and radiation variables, each based on Local Lax-Friedrichs flux. Because we are only focused on time integration and want to test more complex problems that may not have an analytical solution, our error is measured with respect to a discrete reference solution that is computed using extremely small time steps. We point out that results are invariant to the reference integration scheme, as long as a sufficiently small time step is chosen (Lie-Trotter and first-order schemes, however, require extremely small time steps to be used as a reference for second-order schemes). A side result of this testing framework is demonstrating that in the temporal domain, even in the case of nonlinear shocks, we are able to achieve robust second-order convergence to the solution of the ODEs in time that arise post-spatial discretization.

\begin{figure}[!htb]
  \centering
  \begin{subfigure}[b]{0.24\textwidth}
    \includegraphics[width=\textwidth]{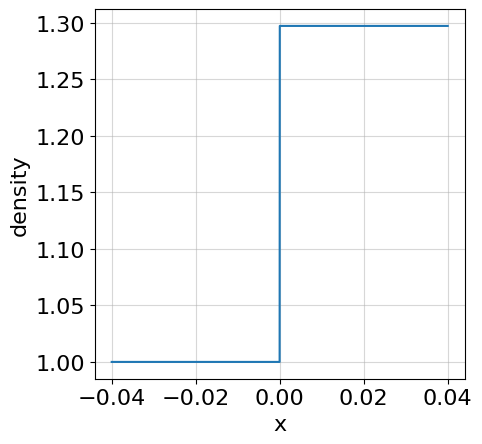}
  \end{subfigure}
  \begin{subfigure}[b]{0.24\textwidth}
    \includegraphics[width=\textwidth]{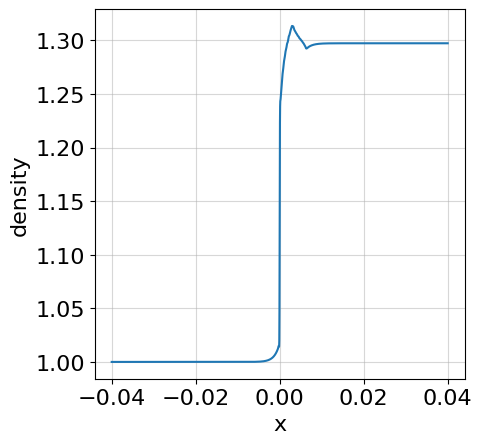}
  \end{subfigure}
  \begin{subfigure}[b]{0.24\textwidth}
    \includegraphics[width=\textwidth]{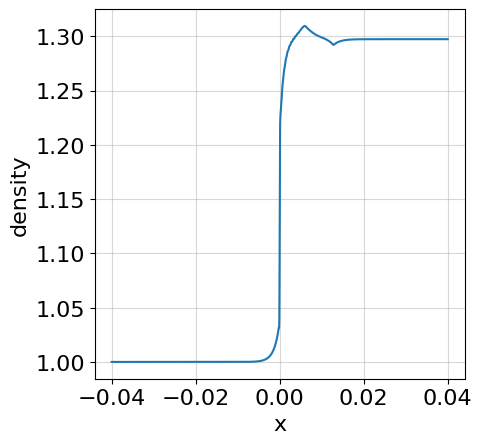}
  \end{subfigure}
  \begin{subfigure}[b]{0.24\textwidth}
    \includegraphics[width=\textwidth]{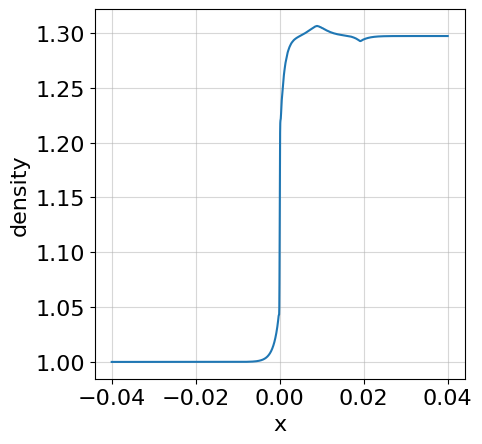}
  \end{subfigure}
  \vspace{-2ex}
  \caption{Mach-1.2 density profiles for $t\in\{0,0.25,0.5,0.75\}$ns}
  \label{fig:m12d}
  \vspace{1ex}
  \begin{subfigure}[b]{0.24\textwidth}
    \includegraphics[width=\textwidth]{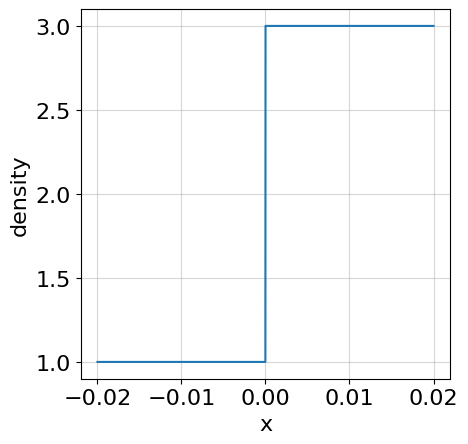}
  \end{subfigure}
  \begin{subfigure}[b]{0.24\textwidth}
    \includegraphics[width=\textwidth]{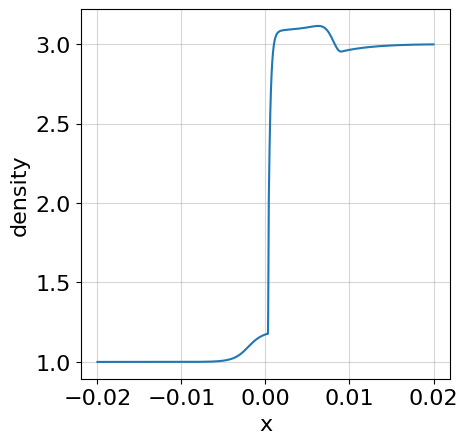}
  \end{subfigure}
  \begin{subfigure}[b]{0.24\textwidth}
    \includegraphics[width=\textwidth]{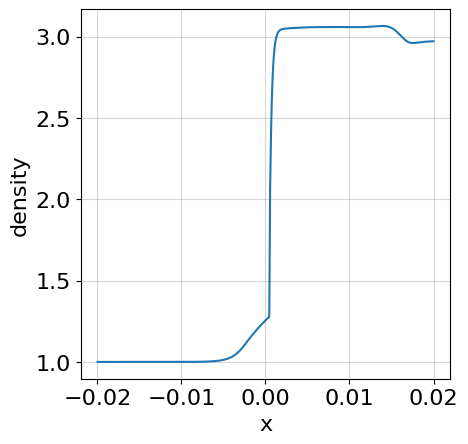}
  \end{subfigure}
  \begin{subfigure}[b]{0.24\textwidth}
    \includegraphics[width=\textwidth]{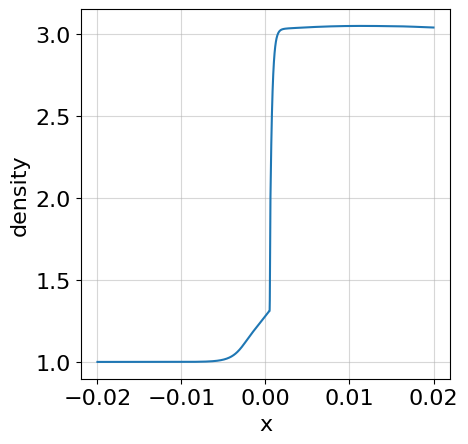}
  \end{subfigure}
  \vspace{-2ex}
  \caption{Mach-3 density profiles for $t\in\{0,0.25,0.5,0.75\}$ns}
  \label{fig:m3d}
  \vspace{1ex}
  \begin{subfigure}[b]{0.24\textwidth}
    \includegraphics[width=\textwidth]{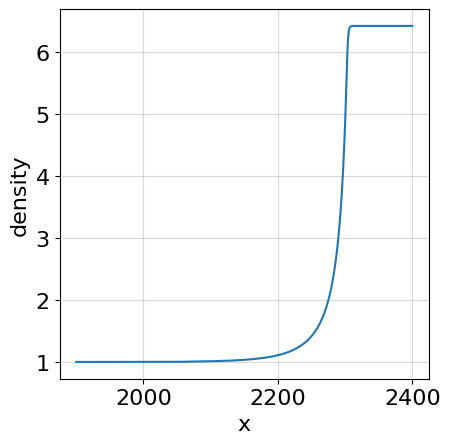}
  \end{subfigure}
  \begin{subfigure}[b]{0.24\textwidth}
    \includegraphics[width=\textwidth]{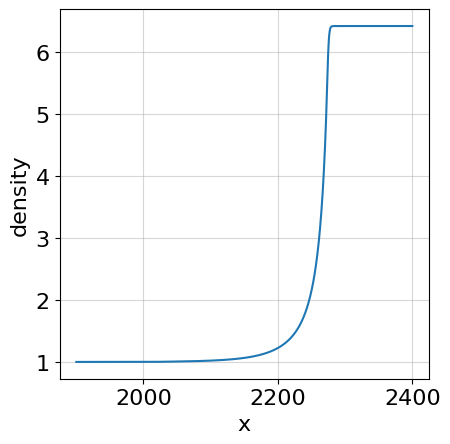}
  \end{subfigure}
  \begin{subfigure}[b]{0.24\textwidth}
    \includegraphics[width=\textwidth]{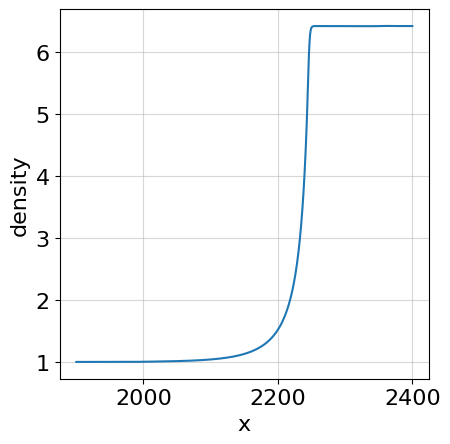}
  \end{subfigure}
  \begin{subfigure}[b]{0.24\textwidth}
    \includegraphics[width=\textwidth]{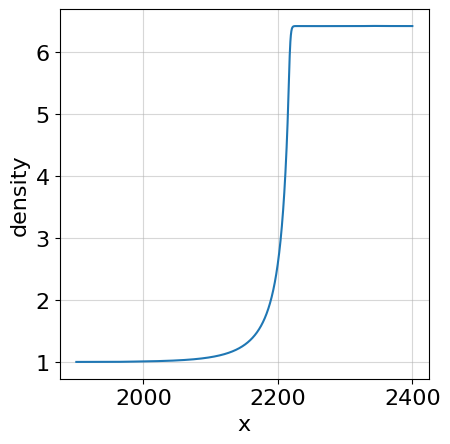}
  \end{subfigure}
  \vspace{-2ex}
  \caption{Mach-45 density profiles for $t\in\{0,50,100,150\}$ns.}
  \label{fig:m45d}
\end{figure}

\begin{figure}[!htb]
  \centering
  \begin{subfigure}[b]{0.24\textwidth}
    \includegraphics[width=\textwidth]{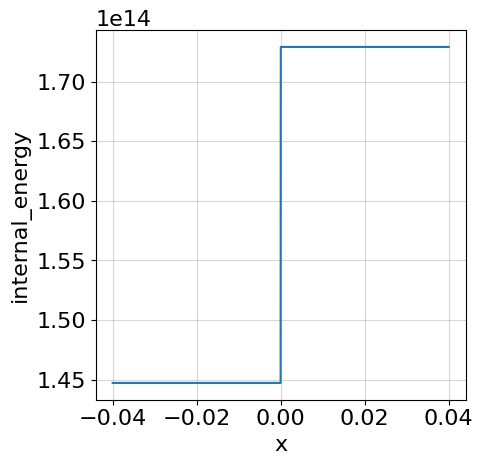}
  \end{subfigure}
  \begin{subfigure}[b]{0.24\textwidth}
    \includegraphics[width=\textwidth]{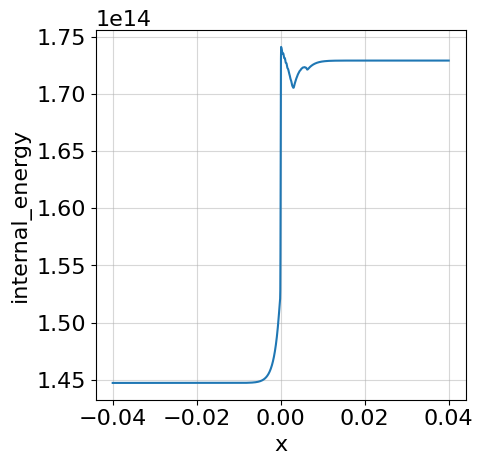}
  \end{subfigure}
  \begin{subfigure}[b]{0.24\textwidth}
    \includegraphics[width=\textwidth]{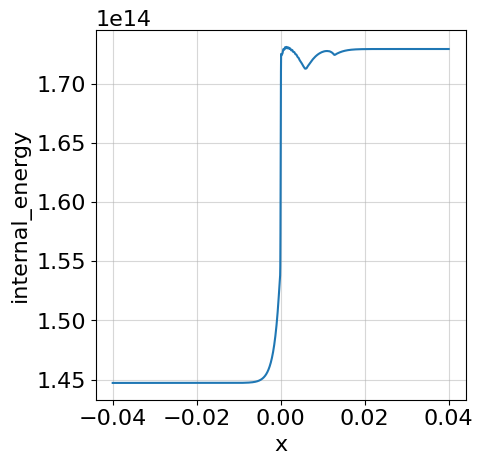}
  \end{subfigure}
  \begin{subfigure}[b]{0.24\textwidth}
    \includegraphics[width=\textwidth]{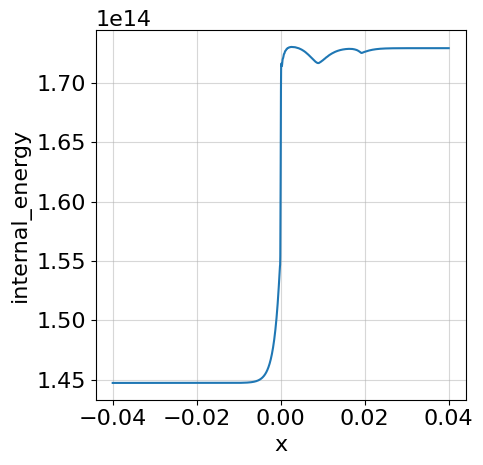}
  \end{subfigure}
  \vspace{-2ex}
  \caption{Mach-1.2 internal energy profiles for $t\in\{0,0.25,0.5,0.75\}$ns}
  \label{fig:m12e}
  \vspace{1ex}
  \begin{subfigure}[b]{0.24\textwidth}
    \includegraphics[width=\textwidth]{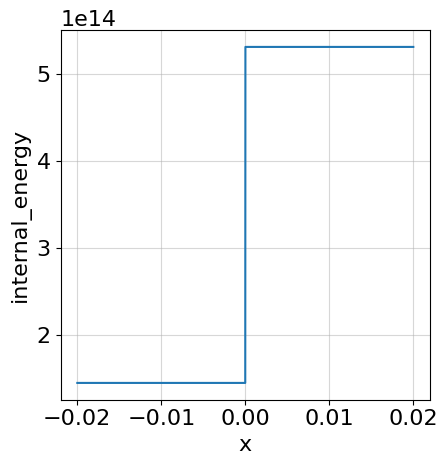}
  \end{subfigure}
  \begin{subfigure}[b]{0.24\textwidth}
    \includegraphics[width=\textwidth]{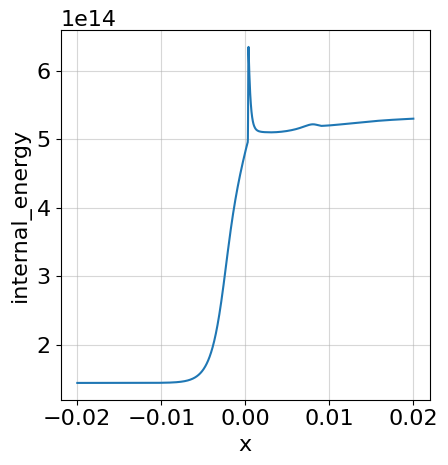}
  \end{subfigure}
  \begin{subfigure}[b]{0.24\textwidth}
    \includegraphics[width=\textwidth]{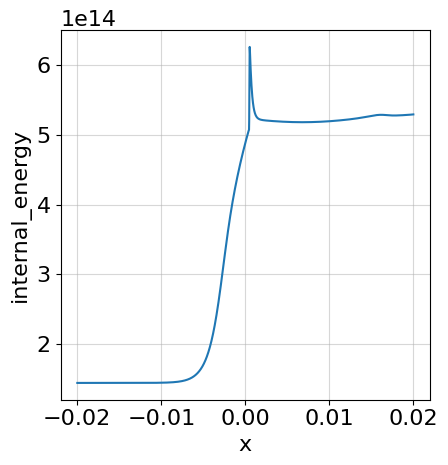}
  \end{subfigure}
  \begin{subfigure}[b]{0.24\textwidth}
    \includegraphics[width=\textwidth]{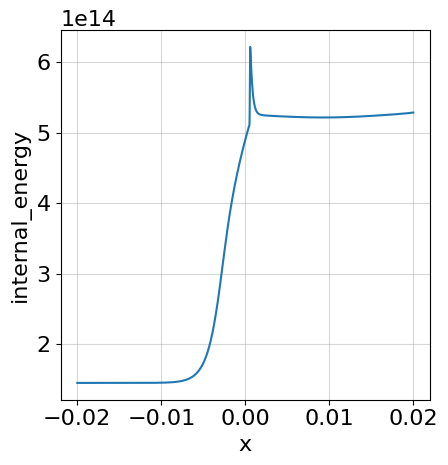}
  \end{subfigure}
  \vspace{-2ex}
  \caption{Mach-3 internal energy profiles for $t\in\{0,0.25,0.5,0.75\}$ns}
  \label{fig:m3e}
  \vspace{1ex}
  \begin{subfigure}[b]{0.24\textwidth}
    \includegraphics[width=\textwidth]{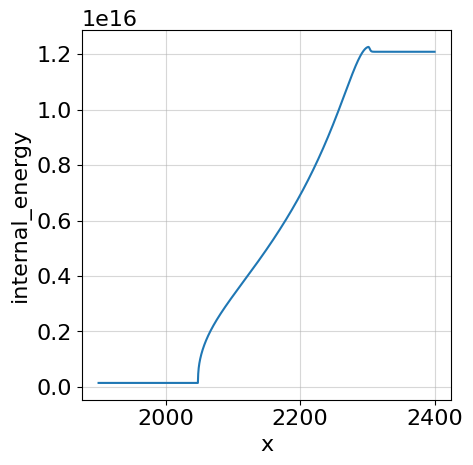}
  \end{subfigure}
  \begin{subfigure}[b]{0.24\textwidth}
    \includegraphics[width=\textwidth]{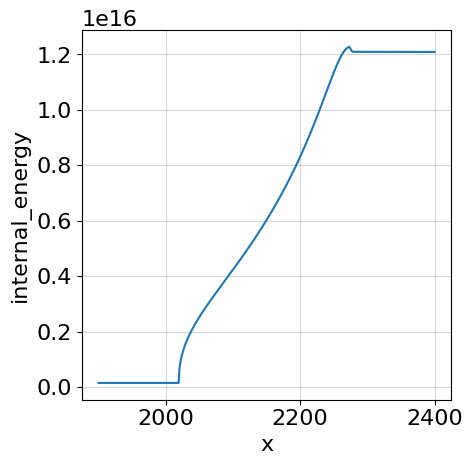}
  \end{subfigure}
  \begin{subfigure}[b]{0.24\textwidth}
    \includegraphics[width=\textwidth]{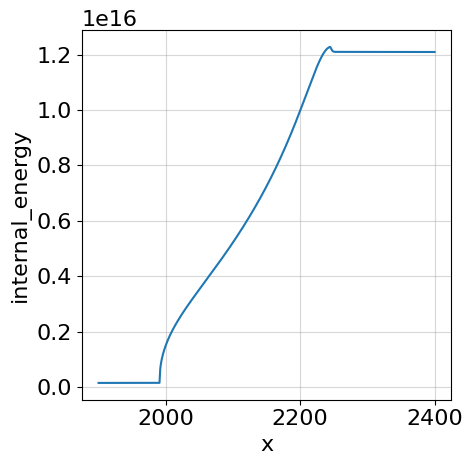}
  \end{subfigure}
  \begin{subfigure}[b]{0.24\textwidth}
    \includegraphics[width=\textwidth]{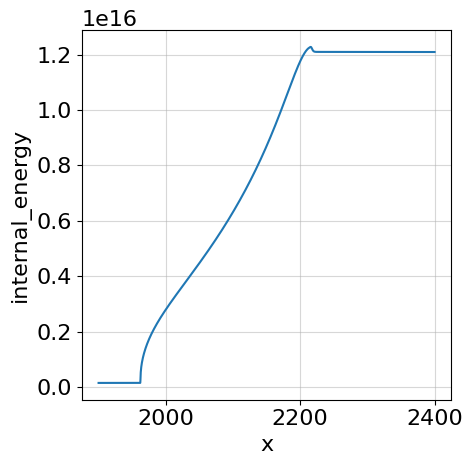}
  \end{subfigure}
  \vspace{-2ex}
  \caption{Mach-45 internal energy profiles for $t\in\{0,50,100,150\}$ns.}
  \label{fig:m45e}
\end{figure}

We test our IMEX schemes using three radiative shock problems, similar to those presented in \cite{Lowrie.2008,Kadioglu.2010,Bolding.2017}. We set up problems with $M=$ 1.2, 3.0, and 45. Table \ref{tab:radshock_specification} shows the material properties:
\begin{center}
\begin{table}[!h]
\centering
\begin{tabular}{ c | c c c}
  & Problem 1& Problem 2 & Problem 3 \\ 
\hline
$M$ & 1.2 &3.0 &45\\
$\gamma$ & 5/3&5/3&5/3 \\
$c_v$[erg/eV/g] &$1.447\times 10^{12}$ &$1.447\times 10^{12}$ &$1.447\times 10^{12}$ \\
$\sigma_a[\mathrm{cm}^{-1}]$ & 577.35 & 577.35 & $4.494\times 10^8 \rho^2/T^{-3.5}$\\
$\sigma_s[\mathrm{cm}^{-1}]$ & 0 & 0 & $0.4006\rho$\\
$\rho_u$[g/cc] & 1.0&1.0&1.0\\
$T_{m,u}$[eV] &100&100&100 \\
$T_{r,u}$[eV] &100&100&100 \\
\end{tabular}
\caption{Radiative shock problem specifications. Subscript ``u'' denotes upstream values.}
\label{tab:radshock_specification}
\end{table}
\end{center}
For problems with $M=1.2$ and $3$ we set shock propagation speed to zero, while the Mach 45 problem has a propagation speed of $-5.7054\times 10^8$ [cm/s] (with respect to the reference frame). We use ideal gas law for EOS. Given $M$, $\rho_u$, $p_u =(\gamma-1)\rho_u c_v T_u$, and $e_u = c_v T_{m,u}$, the downstream values can be evaluated via solving the following Rankine–Hugoniot relations (or jump conditions):
\begin{subequations}\label{eq:jump-full}
\begin{align}
    \rho_u u_u &= \rho_d u_d,\\
    \rho_u u^2_u + p_u + p_{r,u} &=\rho_d u^2_d + p_d + p_{r,d} \\
    (E_u + p_u + p_{r,u})u_u &= (E_d + p_d + p_{r,d})u_d 
\end{align}
\end{subequations}
Generally, downstream values can be evaluated solving the nonlinear system (\ref{eq:jump-full}). However,  when low energy density approximation (e.g., ignoring radiation pressure and radiation energy density in the jump condition) is applied, the following explicit expression can be obtained, 
\begin{subequations}
\begin{align}
    \rho_d &= \rho_u\frac{(\gamma+1)M^2}{2+(\gamma-1)M^2},\\
    u_d &= u_u\frac{\rho_u}{\rho_d},\\
    p_d &= p_u\frac{1+2 \gamma \left(M^2-1\right)}{\gamma+1}.
\end{align}
\end{subequations}
We employ 2nd-order, cell-centered finite volume spatial discretization. Rusanov and central differencing schemes are used for evaluating advective and diffusive numerical fluxes, respectively. Finally, a least-square with Barth-Jespersen limiter is used for high-order flux reconstruction. Reference density and internal energy profiles at four times for each problem are shown in \Cref{fig:m12d,fig:m3d,fig:m45d,fig:m12e,fig:m3e,fig:m45e}. \tcb{All simulation are run with 2048 spatial cells in 1d.}

\tcb{Convergence studies are run starting with timesteps approximately given by the hydrodynamics CFL. If a given integration scheme is not shown on a plot, it means that the scheme in question was unstable at this initial time step size. The Mach-1.2 and Mach 3 problems are integrated to 1 ns for measuring error and convergence, and the Mach-45 problem is integrated to 200 ns. These time intervals are sufficiently long to capture significant interesting and complex dynamics, while also ensuring the solutions have not reached some form of steady or equilibrium solution, where integration error made in early dynamics may no longer be observable. Error is measured as a discrete relative $\ell^2$-, $\ell^1$-, and $\ell^\infty$-error of the solution over the physical domain with respect to a reference solution computed using an extremely small time step. We use multiple norms to demonstrate that the convergence is robust across different measures, even an $\infty$-norm, a common choice for shock problems.  We use SSP-LDIRK3(3,3,2) (see \Cref{sec:results:imex}) to compute the reference solutions as it is consistently the most accurate integration scheme, and we use a timestep $10\times$ smaller than the smallest one shown in convergence plots to ensure the reference error is $\ll$ than the solution being compared against.}


\subsection{IMEX schemes}\label{sec:results:imex}

{\color{black}
For IMEX schemes, we consider the Lie-Trotter operator split with classical semi-implicit Euler for radiation and a forward Euler (\emph{Op-Split}) or a 3rd-order explicit total variation diminishing RK method (\emph{Op-Split-TVD3}) for hydrodynamics. For LIMEX schemes, we have tested many from the literature, e.g. \cite{Kennedy.2003tv4,Ascher.1997,Pareschi.2005,Giraldo.2013,Conde.2017,Boscarino.2016,Sebastiano.2023}. In addition to the simplest semi-implicit first-order scheme \eqref{eq:limex-o1} we consider four additional schemes, three second order and one third order, that have consistently performed superior to others, particularly in terms of stability. Avoiding instability and nonphysical negativities has proved not trivial for higher-order integration schemes. All of of the schemes we have chosen have the same number of implicit and explicit stages, and thus do \emph{not} have an ESDIRK \cite{Kennedy.2003tv4} or padded tableaux structure \cite{Ascher.1997}. In fact, there is theoretical precedence for this -- in \cite{Pareschi.2005} it is proven that for IMEX schemes to be asymptotic preserving for hyperbolic equations with stiff relaxation terms, they must have equal implicit and explicit number of stages/not have explicit stages in the implicit method. Radiation hydrodynamics \eqref{eq:rad-hydro} can indeed be seen as a more general form of a hyperbolic equation with relaxation -- we couple the hyperbolic hydrodynamics equations with stiff radiation energy and temperature equations, where the stiff relaxation parameter (typically referred to as $1/\epsilon$) corresponds to the radiation and temperature coupling coefficients, $\sigma_E c$ and $\sigma_p a c$. Three of the IMEX schemes we choose were developed in \cite{Pareschi.2005} specifically to be asymptotic preserving for hyperbolic equations with stiff relaxation terms. Intuitively, these schemes do not allow explicit evaluations of the implicit operator, which, even if damped by a following implicit stage, seem more likely to induce negativities or instabilities. We also note that most third-order schemes we tested have had poor stability properties for our problems, and we report on the only scheme we have tested that has proven robust. We believe this is due to a higher order method also having higher probability for nonphysical solution values such as negative temperatures, in problems with shocks and sharp gradients; although the code ``fixes'' any negativities to avoid immediate breakdown, too much of this behavior can lead to poor accuracy or numerical instabilities. 

\textbf{H-LDIRK2(2,2,2)} \cite[Table II]{Pareschi.2005} (2nd-order scheme with 2nd-order L-stable DIRK implicit method and 2nd-order SSP explicit):
\begin{align*}
    \begin{array}{c | c c c}
    0 & 0 & 0 \\
    1 & 1 & 0 \\\hline
    & 1/2 & 1/2
    \end{array},
\hspace{5ex}
    \begin{array}{c | c c}
    \gamma & \gamma & 0 \\
    1 -\gamma & 1-2\gamma & \gamma \\\hline
    & 1/2 & 1/2
    \end{array},
\hspace{2ex} \textnormal{for }\gamma=1-1/\sqrt{2}.
\end{align*}
\noindent\textbf{SSP-LDIRK2(3,3,2)} \cite[Table IV]{Pareschi.2005} (2nd-order scheme with 2nd-order L-stable DIRK implicit method and 2nd-order SSP explicit):
\begin{align*}
    \begin{array}{c | c c c}
    0 & 0 & 0 & 0 \\
    1/2 & 1/2 & 0 & 0\\
    1 & 1/2 & 1/2 & 0 \\\hline
    & 1/3 & 1/3 & 1/3 
    \end{array},
\hspace{5ex}
    \begin{array}{c | c c c}
    0.25 & 0.25 & 0 & 0 \\
    0.25 & 0 & 0.25 & 0 \\
    1 & 1/3 & 1/3 & 1/3 \\\hline
    & 1/3 & 1/3 & 1/3 
    \end{array}.
\end{align*}
\textbf{SSP-LDIRK3(3,3,2)} \cite[Table V]{Pareschi.2005} (2nd-order scheme with 2nd-order L-stable DIRK implicit method and 3rd-order SSP explicit):
\begin{align*}
    \begin{array}{c | c c c}
    0 & 0 & 0 & 0 \\
    1 & 1 & 0 & 0\\
    1/2 & 1/4 & 1/4 & 0 \\\hline
    & 1/6 & 1/6 & 2/3
    \end{array},
\hspace{5ex}
    \begin{array}{c | c c c}
    \gamma & \gamma & 0 & 0 \\
    1-\gamma & 1-2\gamma & \gamma & 0 \\
    1/2 & 1/2 - \gamma & 0 & \gamma \\\hline
    & 1/6 & 1/6 & 2/3
    \end{array},
    \hspace{2ex} \textnormal{for }\gamma=1-1/\sqrt{2}.
\end{align*}

\textbf{I-IMEX(3,4,3)} \cite[Eq. (30)]{Sebastiano.2023} (3rd-order scheme):
\begin{align*}
    &\begin{array}{c | c c c c}
        0 & 0 & 0 & 0 & 0 \\
        0.4358665215 & 0.4358665215 & 0 & 0 & 0 \\
        0.7179332608 & 1.243893189 & -0.5259599287 & 0 & 0 \\
        1 & 0.6304125582 & 0.7865807402 & -0.4169932983 & 0 \\\hline
        & 0 & 1.208496649 & -0.644363171 & 0.4358665215 \\
    \end{array}, \\
    &\begin{array}{c | c c c c}
        0.4358665215 & 0.4358665215 & 0 & 0 & 0 \\
        0.4358665215 & 0 & 0.4358665215 & 0 & 0 \\
        0.7179332608 & 0 & 0.2820667392 & 0.4358665215 & 0 \\
        1 & 0 & 1.208496649 & -0.644363171 & 0.43586652150 \\\hline
         & 0 & 1.208496649 & -0.644363171 & 0.4358665215
    \end{array},
\end{align*}
}

\subsection{Fine spatial mesh}
\begin{figure}[!htb]
  \centering
  \begin{subfigure}[b]{0.7\textwidth}
    \includegraphics[width=\textwidth]{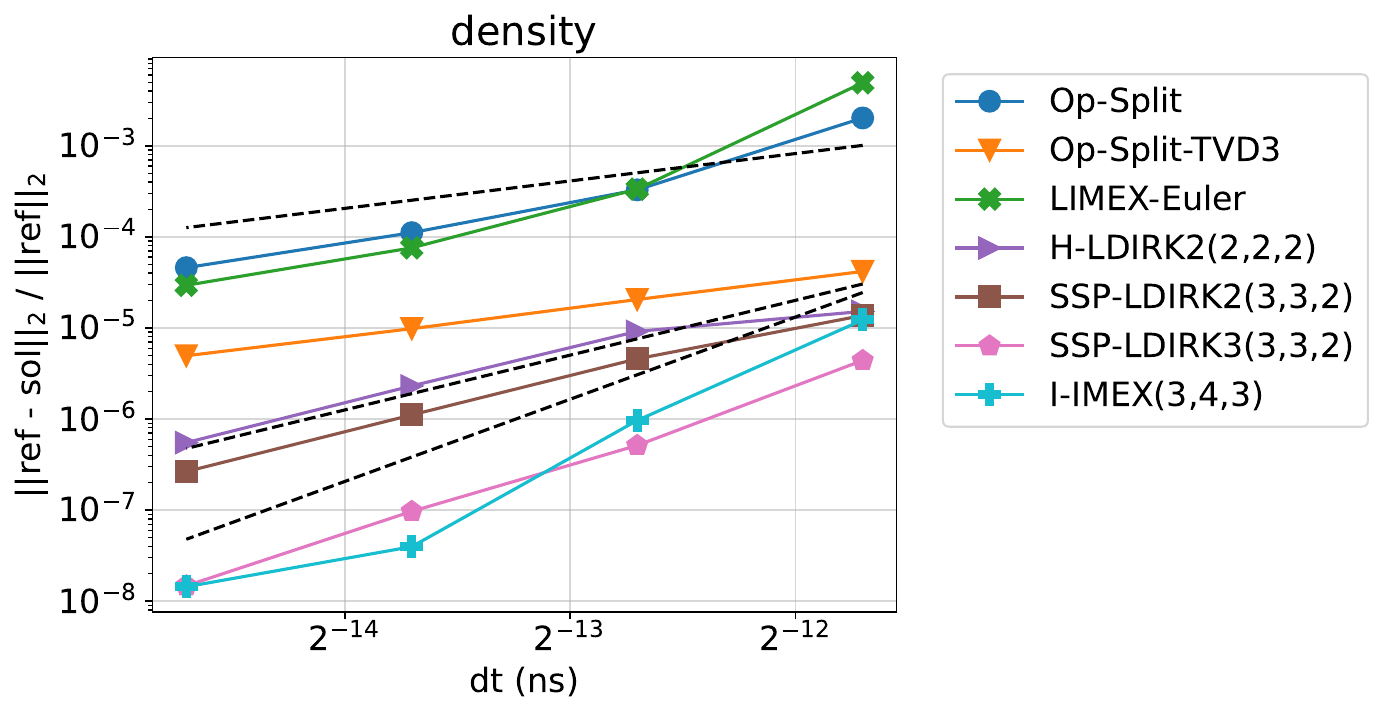}
  \end{subfigure}
  \\
  \begin{subfigure}[b]{0.45\textwidth}
    \includegraphics[width=\textwidth]{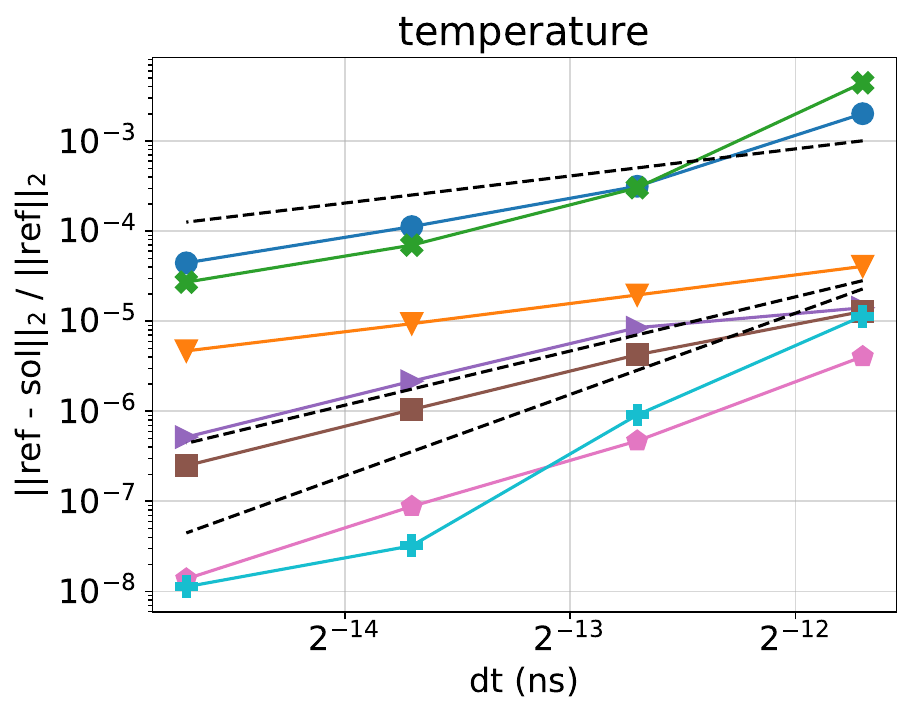}
  \end{subfigure}
  \begin{subfigure}[b]{0.45\textwidth}
    \includegraphics[width=\textwidth]{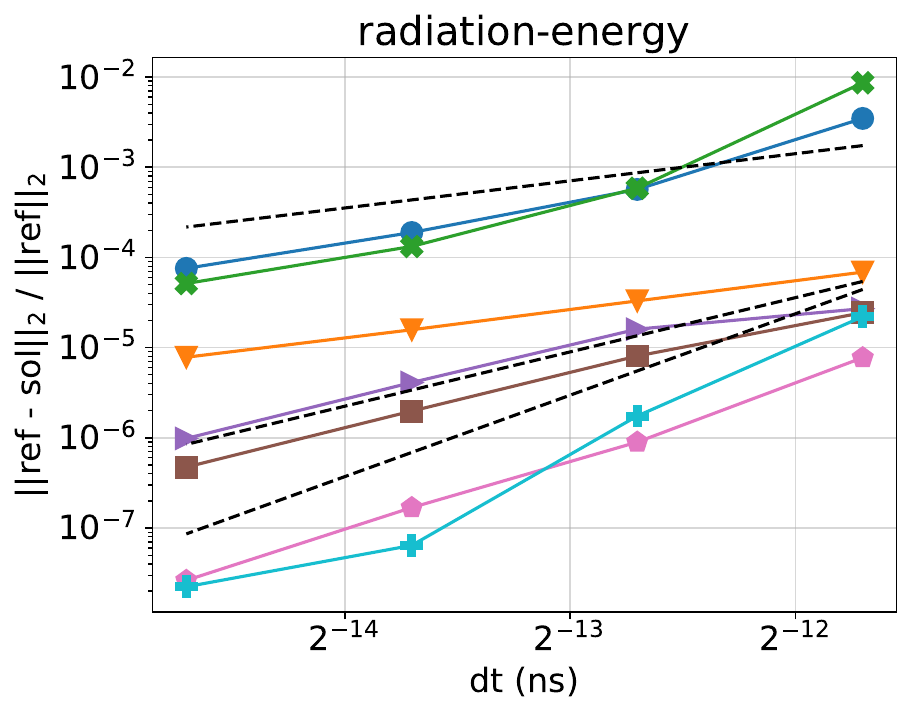}
  \end{subfigure}
    \caption{Mach-1.2 relative $\ell^2$-error; dotted lines demonstrate first, second, and third order convergence.}
      \label{fig:m1p2}
\end{figure}
\begin{figure}[!htb]
  \centering
  \begin{subfigure}[b]{0.7\textwidth}
    \includegraphics[width=\textwidth]{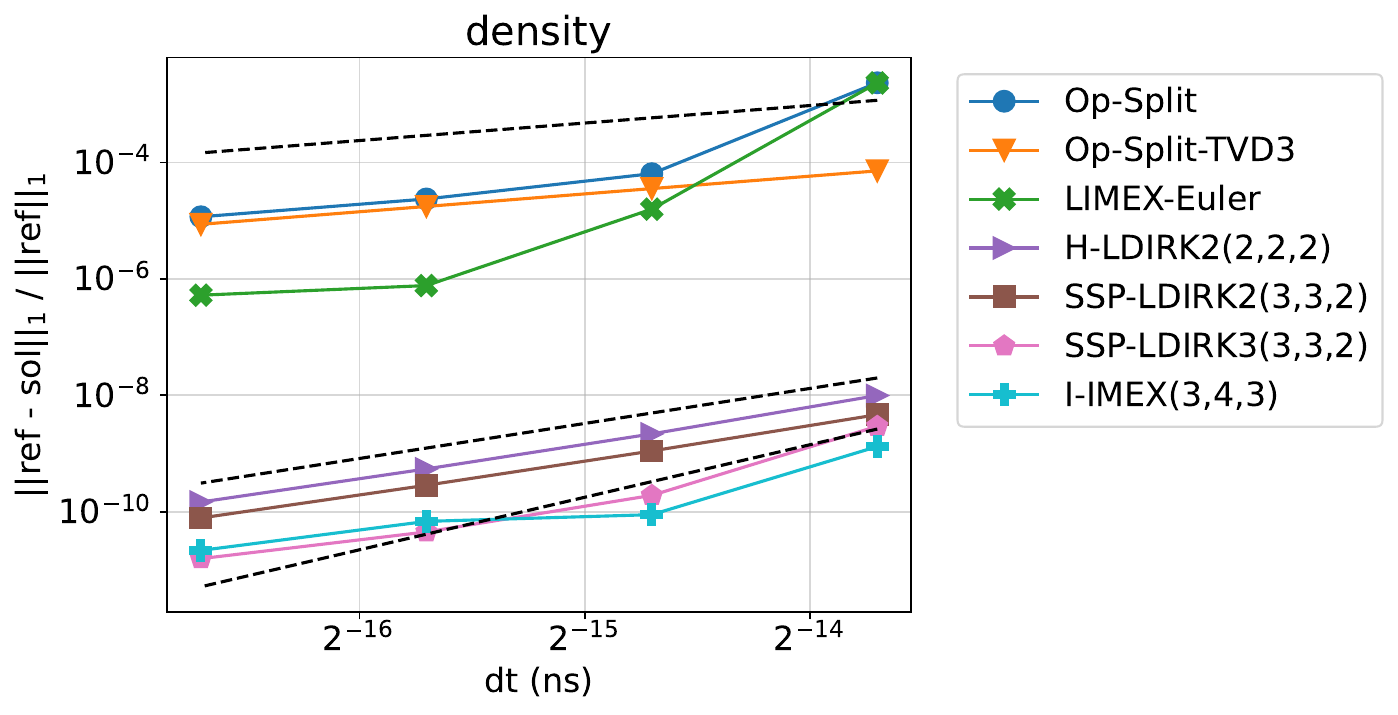}
  \end{subfigure}
  \\
  \begin{subfigure}[b]{0.45\textwidth}
    \includegraphics[width=\textwidth]{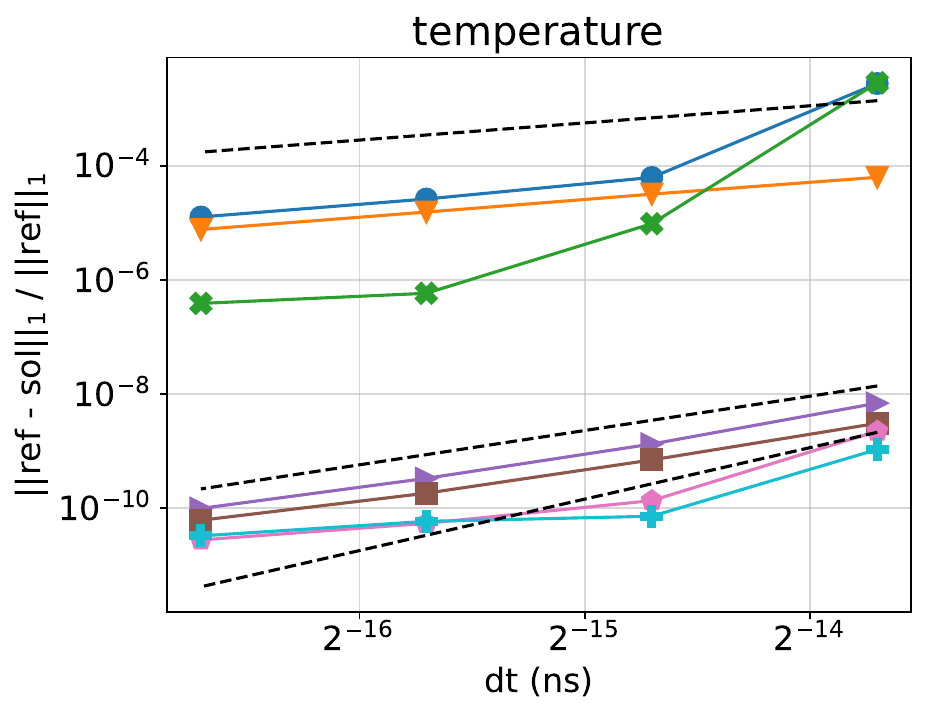}
  \end{subfigure}
  \begin{subfigure}[b]{0.45\textwidth}
    \includegraphics[width=\textwidth]{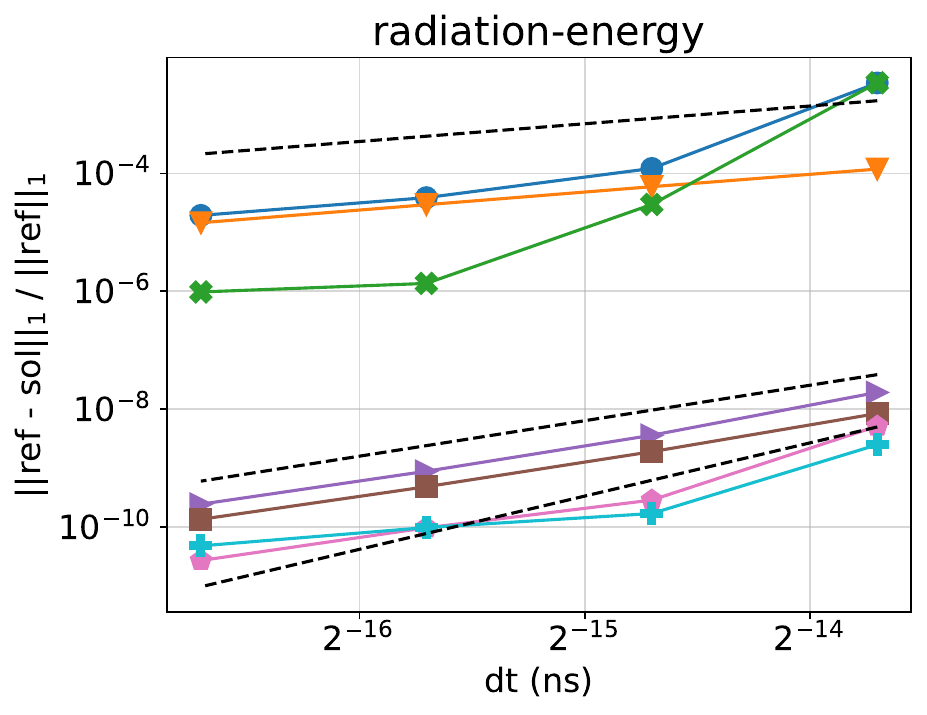}
  \end{subfigure}
    \caption{Mach-3 relative $\ell^1$-error; dotted lines demonstrate first, second, and third order convergence.}
      \label{fig:m3}
\end{figure}
\begin{figure}[!htb]
  \centering
  \begin{subfigure}[b]{0.7\textwidth}
    \includegraphics[width=\textwidth]{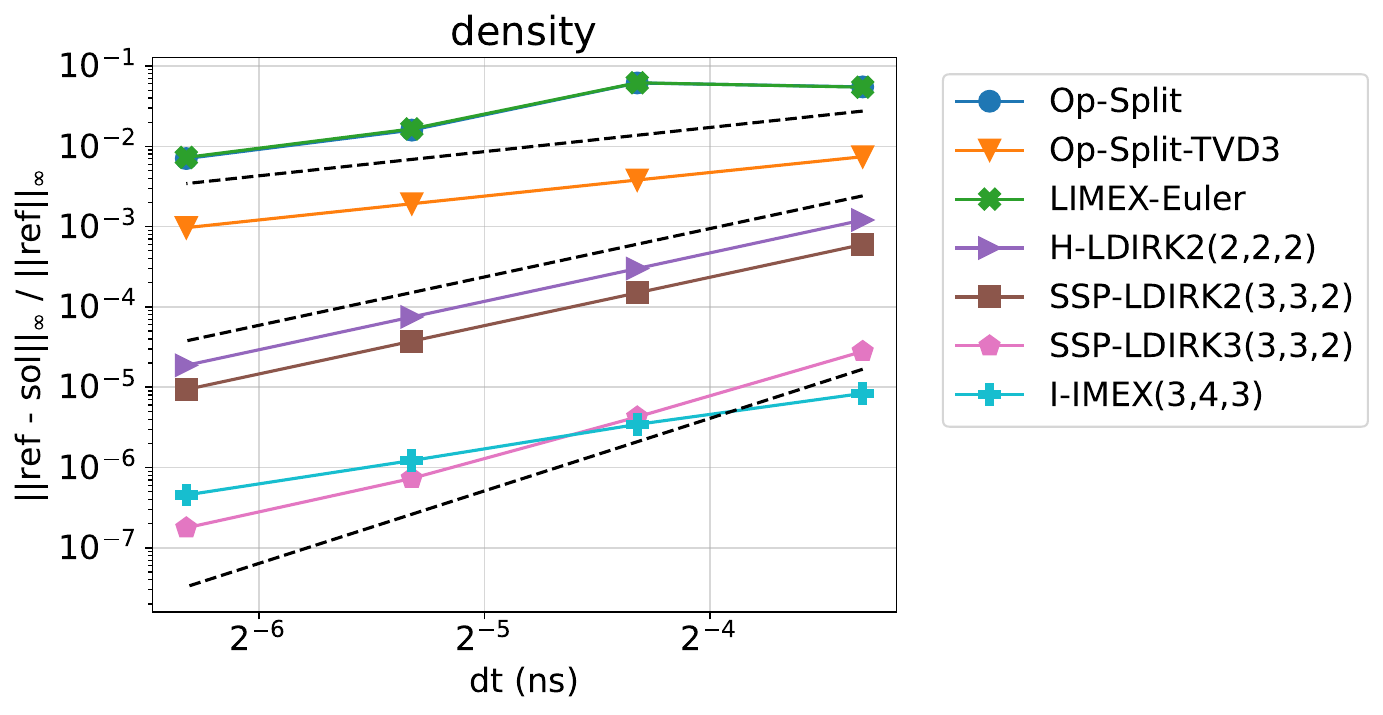}
  \end{subfigure}
  \\
  \begin{subfigure}[b]{0.45\textwidth}
    \includegraphics[width=\textwidth]{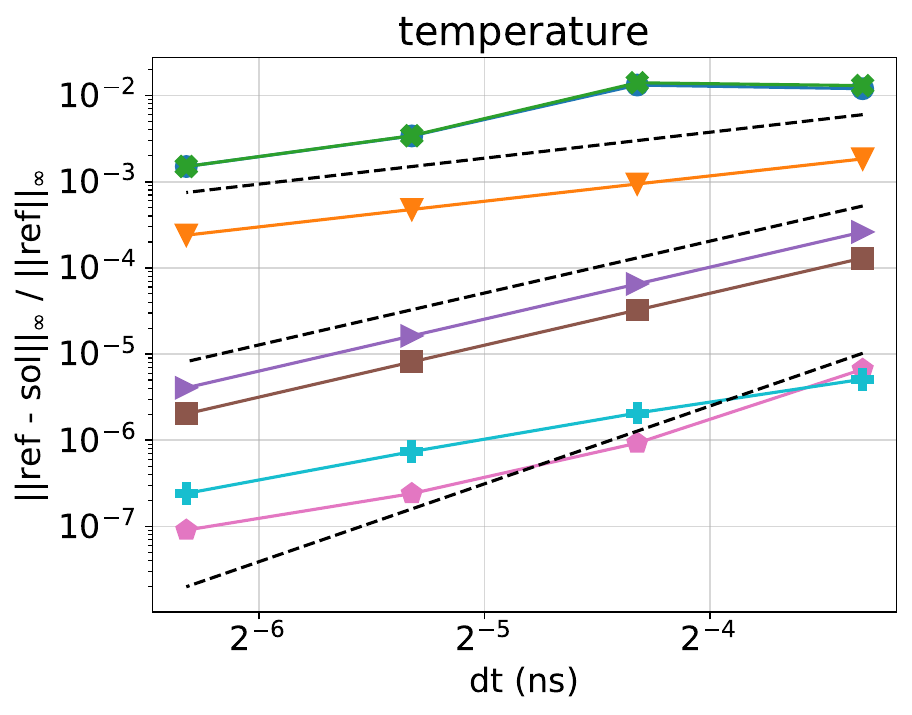}
  \end{subfigure}
  \begin{subfigure}[b]{0.45\textwidth}
    \includegraphics[width=\textwidth]{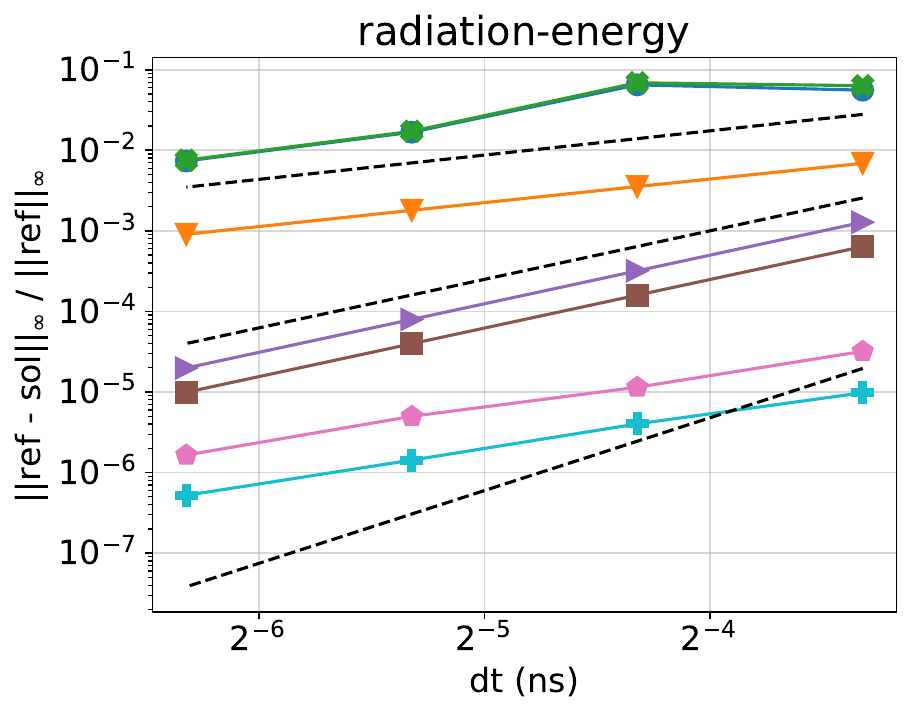}
  \end{subfigure}
    \caption{Mach-45 relative $\ell^\infty$ error; dotted lines demonstrate first, second, and third order convergence.}
      \label{fig:m45}
\end{figure}

{\color{black}
We begin by presenting results on a fine spatial mesh of 2000 cells in 1d for all three problems in \Cref{fig:m1p2,fig:m3,fig:m45}. We choose the density, radiation energy, and temperature as representative variables in hydrodynamics, radiation, and their coupling, respectively, but note that similar behavior is observed for momentum and hydro total energy. We also present one problem in relative $\ell^1$- (Mach-3, \Cref{fig:m3}), $\ell^2-$ (Mach-1.2, \Cref{fig:m1p2}), and $\ell^\infty$- (Mach-45, \Cref{fig:m45}) error; as before, analogous results hold for all problems in all three different norms, and we present this way to demonstrate our results in multiple norms.

The most important observation in our results is that for large time steps, the proposed IMEX schemes achieve accuracy in all variables orders of magnitude smaller than the operator split methods, even when a 3rd-order explicit SSP method is used within the operator split. This is particularly pronounced for the Mach-3 problem (\Cref{fig:m3}), where IMEX schemes achieve \emph{4--5 orders of magnitude} smaller error than the operator split methods. In almost all cases we are also able to observe second order accuracy with refinement in time step. For Mach-1.2, we observe third order convergence for I-IMEX(3,4,3) as well as SSP-LDIRK3(3,3,2). The latter scheme is only formally second order, but it appears that we capture second order differentials in radiation and coupling with sufficiently high accuracy that our error is dominated by hydrodynamics, and SSP-LDIRK3(3,3,2) is built on a 3rd-order explicit integrator. In Mach-3 (\Cref{fig:m3}) and Mach-45 (\Cref{fig:m45}), each of these schemes demonstrate order reduction, but we believe this is due to implicit solve tolerance in the code rather than the actual methods (as well as due to coefficient accuracy for I-IMEX(3,4,3), which is only provided to 10 digits in \cite{Sebastiano.2023}); i.e., the relative error is already very small at the largest time steps, and inexact linear and nonlinear implicit solves in the implicit radiation equations will limit how much more accurate these methods can be. 

Finally, we compare the wallclock time of the different methods for all three problems. To do this, we measure the wallclock time of each integration scheme over the four time steps considered in \Cref{fig:m1p2,fig:m3,fig:m45} and plot it with respect to the wallclock time of LIMEX-Euler (see \Cref{fig:wallclock}). We choose LIMEX Euler because our IMEX code is not fully optimized to avoid repeat function calls or unnecessary stage storage, so the operator split is $\approx 25\%$ faster, despite being more or less identical on paper. Aside from this optimization, in most cases we see perfect linear scaling in cost with respect to stage, that is, the 2-stage scheme takes $2\times$ longer than LIMEX-Euler, and analogously for 3-stage and 4-stage. Mach-45 appears to offer sublinear scaling in some cases, which appears to be due to the implicit stage equations with more stages being slightly easier to solve. Comparing with the error plots in \Cref{fig:m1p2,fig:m3,fig:m45}, we see that for a $2-4\times$ increase in computational cost, one can achieve hundreds or thousands times better accuracy.

\begin{figure}[!htb]
    \centering
    \includegraphics[width=\textwidth]{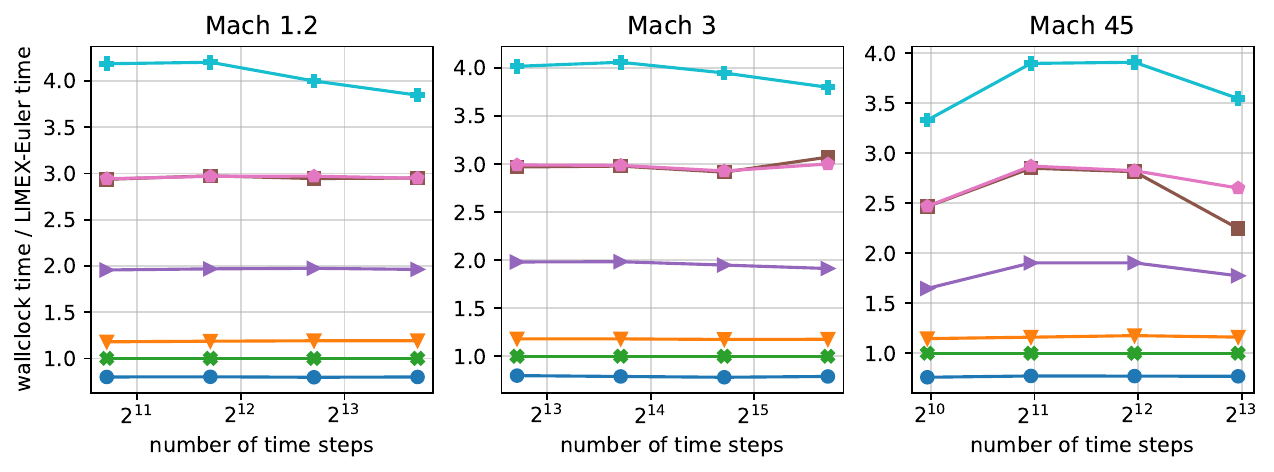}
    \caption{Wallclock time for simulations in \Cref{fig:m1p2,fig:m3,fig:m45}.}
      \label{fig:wallclock}
\end{figure}

}

\subsection{Coarse spatial mesh}

{\color{black}
For our last test, we take the same problems on a coarse spatial mesh of 200 cells in 1d, closer to resolution that could be afforded in 2d or 3d calculations. Now we have shocks and gradients that are definitely not resolved by the spatial mesh. We run similar convergence studies as before, starting with a largest timestep defined by the approximate coarse-mesh CFL, and refining by a factor of two. Results are shown in \Cref{fig:coarse}. Remarkably, we are still able to observe excellent convergence, and actually achieve more robust 3rd-order accuracy with I-IMEX(3,4,3) (in part, this may be because the relative error over timesteps considered is larger, so we do not run into issues with implicit solve tolerance or RK coefficient accuracy). However, we emphasize that here we are measuring accuracy and convergence with respect to the spatially discretized problem, not the underlying PDE solution, so we are not claiming 3rd order accuracy of shocks with respect to physical solution.

\begin{figure}[!htb]
  \centering
      \begin{subfigure}[b]{0.45\textwidth}
    \includegraphics[width=\textwidth]{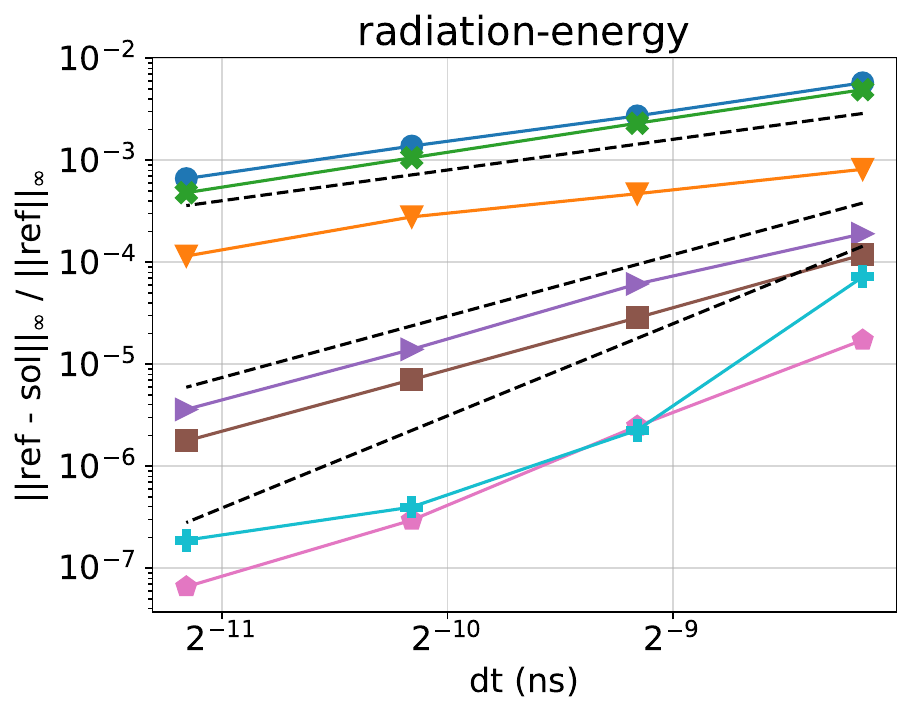}
    \caption{Mach-1.2 relative $\ell^\infty$ error.}
  \end{subfigure}
  \begin{subfigure}[b]{0.45\textwidth}
    \includegraphics[width=\textwidth]{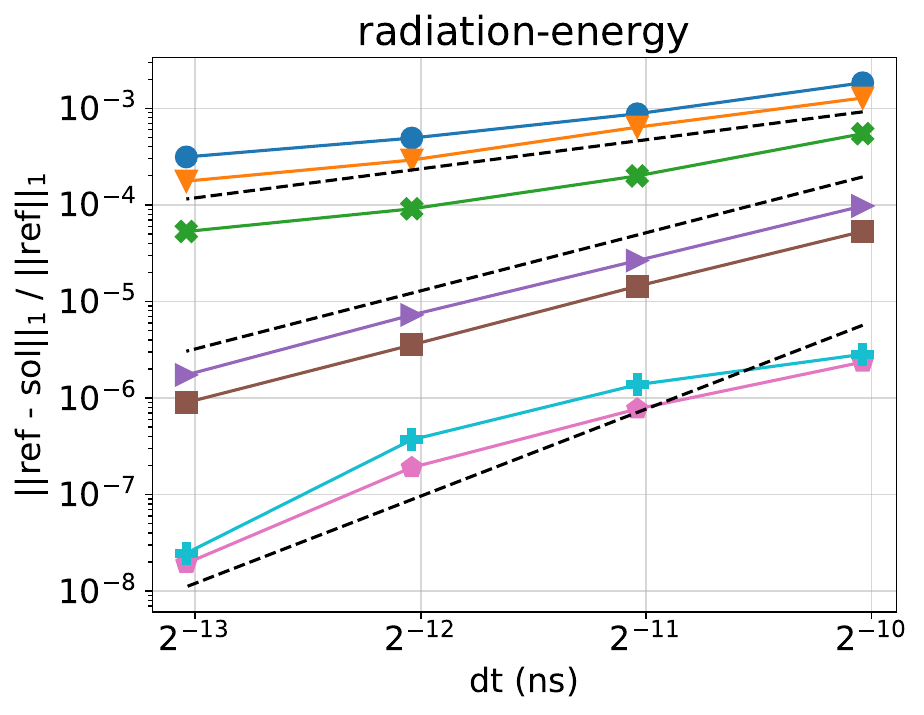}
    \caption{Mach-3 relative $\ell^1$ error.}
  \end{subfigure}
  \begin{subfigure}[b]{0.45\textwidth}
    \includegraphics[width=\textwidth]{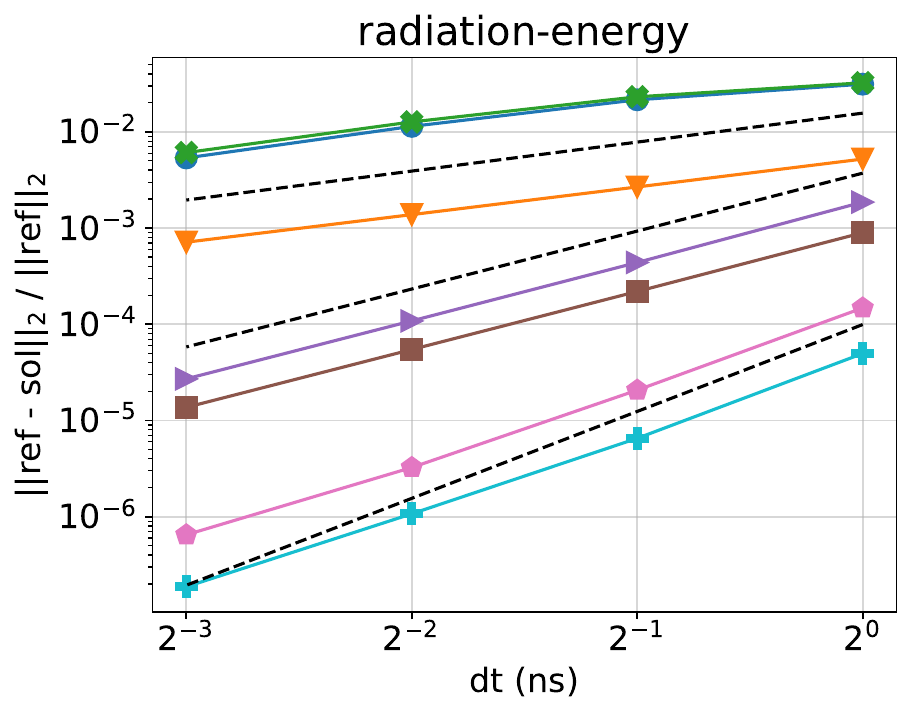}
    \caption{Mach-45 relative $\ell^2$ error.}
  \end{subfigure}
  \caption{Radiation energy relative $\ell^\infty,\ell^1$- and $\ell^2$-error in the three problems for coarse spatial mesh.}
  \label{fig:coarse}
\end{figure}

}

\section{Conclusions}\label{sec:conc}

We introduce a general implicit-explicit time integration framework for radiation hydrodynamics with grey radiation diffusion. The framework is based on the LIMEX methods from \cite{Boscarino.2015}, and provides a natural framework to obtain second or higher-order accuracy in time, without iterating between the hydrodynamics and radiation equations as in \cite{Kadioglu.2010,Bolding.2017}. We prove conservation of energy for the integration schemes, and on 1d radiative shock problems demonstrate 2nd-order, and in some cases 3rd-order, convergence in radiation and hydrodynamics variables. In addition the error constants for timesteps on the hydrodynamics timescales are orders of magnitude lower than the error obtained through a classical Lie-Trotter operator split.

\tcb{Ongoing work will extend the framework to a full kinetic thermal radiation transport (TRT) model, coupling with our recent work on IMEX for TRT \cite{imex-trt}, and also consider integrating in a Lagrangian reference frame. Although we see no conceptual difficulties extending our framework, it is possible the full TRT coupling will provide additional stiffness and accuracy challenges compared with the simpler diffusion model considered here, which can be difficult for higher-order integrators, something we have observed with IMEX for TRT \cite{imex-trt}. We are also exploring adaptivity in time, one of the natural benefits that can be facilitated by holistic RK-type integrators, e.g. \cite{Soderlind.2006,Ranocha.2022}, in comparison with a first-order operator split, and utilizing a new specialized class of integrators called nonlinearly partitioned Runge--Kutta (NPRK) methods \cite{nprk1,nprk2} we have developed specifically for nonlinear partitioning of multiphysics, such as we have considered here. A major challenge of incorporating additional physics is that it may require additional ($>2$) partitions or scale-bridging, and we believe that NPRK methods will provide a clean framework in which to develop such integrators. Finally, we emphasize that the observed accuracy and convergence measured in this paper, as well as any adaptivity in time, is with respect to the spatially discretized problem, not the physical PDE solution. Long-term, we hope to couple the spatial and temporal discretizations for holistic, high-order, adaptive simulations coupling radiation and hydrodynamics, and potentially additional models.}

\section*{Acknowledgements}
This work was supported by the Laboratory Directed Research and Development program of Los Alamos National Laboratory under project number 20220174ER. The authors thank Steven Roberts from Lawrence Livermore National Laboratory for his helpful guidance on three-way ARK methods discussed in \Cref{sec:time:ark}. This research used resources provided by the Darwin testbed at Los Alamos National Laboratory (LANL) which is funded by the Computational Systems and Software Environments subprogram of LANL's Advanced Simulation and Computing program (NNSA/DOE). LA-UR-23-24702.

\bibliographystyle{elsarticle-num} 
\bibliography{refs.bib}

\begin{thebibliography}{10}
\expandafter\ifx\csname url\endcsname\relax
  \def\url#1{\texttt{#1}}\fi
\expandafter\ifx\csname urlprefix\endcsname\relax\def\urlprefix{URL }\fi
\expandafter\ifx\csname href\endcsname\relax
  \def\href#1#2{#2} \def\path#1{#1}\fi

\bibitem{Kadioglu.2010}
S.~Y. Kadioglu, D.~A. Knoll, R.~B. Lowrie, R.~M. Rauenzahn, {A second order
  self-consistent IMEX method for radiation hydrodynamics}, Journal of
  Computational Physics 229~(22) (2010) 8313--8332.
\newblock \href {https://doi.org/10.1016/j.jcp.2010.07.019}
  {\path{doi:10.1016/j.jcp.2010.07.019}}.

\bibitem{Bolding.2017}
S.~Bolding, J.~Hansel, J.~D. Edwards, J.~E. Morel, R.~B. Lowrie, {Second-order
  discretization in space and time for radiation-hydrodynamics}, Journal of
  Computational Physics 338 (2017) 511--526.
\newblock \href {https://doi.org/10.1016/j.jcp.2017.02.063}
  {\path{doi:10.1016/j.jcp.2017.02.063}}.

\bibitem{strang1968construction}
G.~Strang, On the construction and comparison of difference schemes, SIAM
  journal on numerical analysis 5~(3) (1968) 506--517.

\bibitem{Haines.2022}
B.~M. Haines, D.~E. Keller, K.~P. Long, M.~D. McKay, Z.~J. Medin, H.~Park,
  R.~M. Rauenzahn, H.~A. Scott, K.~S. Anderson, T.~J.~B. Collins, L.~M. Green,
  J.~A. Marozas, P.~W. McKenty, J.~H. Peterson, E.~L. Vold, C.~D. Stefano,
  R.~S. Lester, J.~P. Sauppe, D.~J. Stark, J.~Velechovsky, {The development of
  a high-resolution Eulerian radiation-hydrodynamics simulation capability for
  laser-driven Hohlraums}, Physics of Plasmas 29~(8) (2022) 083901.
\newblock \href {https://doi.org/10.1063/5.0100985}
  {\path{doi:10.1063/5.0100985}}.

\bibitem{fuksman2021two}
J.~D.~M. Fuksman, H.~Klahr, M.~Flock, A.~Mignone, A two-moment radiation
  hydrodynamics scheme applicable to simulations of planet formation in
  circumstellar disks, The Astrophysical Journal 906~(2) (2021) 78.

\bibitem{einkemmer2018efficient}
L.~Einkemmer, M.~Moccaldi, A.~Ostermann, Efficient boundary corrected strang
  splitting, Applied Mathematics and Computation 332 (2018) 76--89.

\bibitem{sportisse2000analysis}
B.~Sportisse, An analysis of operator splitting techniques in the stiff case,
  Journal of computational physics 161~(1) (2000) 140--168.

\bibitem{zingale2019improved}
M.~Zingale, M.~Katz, J.~Bell, M.~Minion, A.~Nonaka, W.~Zhang, Improved coupling
  of hydrodynamics and nuclear reactions via spectral deferred corrections, The
  Astrophysical Journal 886~(2) (2019) 105.

\bibitem{zingale2022improved}
M.~Zingale, M.~Katz, A.~Nonaka, M.~Rasmussen, An improved method for coupling
  hydrodynamics with astrophysical reaction networks, The Astrophysical Journal
  936~(1) (2022) 6.

\bibitem{Rice.1960}
J.~R. Rice, {Split Runge-Kutta method for simultaneous equations}, Journal of
  Research of the National Bureau of Standards Section B Mathematics and
  Mathematical Physics 64B~(3) (1960) 151.
\newblock \href {https://doi.org/10.6028/jres.064b.018}
  {\path{doi:10.6028/jres.064b.018}}.

\bibitem{hofer1976partially}
E.~Hofer, A partially implicit method for large stiff systems of odes with only
  few equations introducing small time-constants, SIAM Journal on Numerical
  Analysis 13~(5) (1976) 645--663.

\bibitem{hairer1981order}
E.~Hairer, Order conditions for numerical methods for partitioned ordinary
  differential equations, Numerische Mathematik 36~(4) (1981) 431--445.

\bibitem{Kennedy.2003tv4}
C.~A. Kennedy, M.~H. Carpenter, {Additive Runge–Kutta schemes for
  convection–diffusion–reaction equations}, Applied Numerical Mathematics
  44~(1-2) (2003) 139--181.
\newblock \href {https://doi.org/10.1016/s0168-9274(02)00138-1}
  {\path{doi:10.1016/s0168-9274(02)00138-1}}.

\bibitem{Hairer.1996ci}
E.~Hairer, G.~Wanner, {Solving Ordinary Differential Equations II, Stiff and
  Differential-Algebraic Problems}, Springer Series in Computational
  Mathematics, Springer Berlin Heidelberg, 1996.
\newblock \href {https://doi.org/10.1007/978-3-642-05221-7}
  {\path{doi:10.1007/978-3-642-05221-7}}.

\bibitem{Fuksman.2019}
J.~D.~M. Fuksman, A.~Mignone, {A Radiative Transfer Module for Relativistic
  Magnetohydrodynamics in the PLUTO Code}, The Astrophysical Journal Supplement
  Series 242~(2) (2019) 20.
\newblock \href {http://arxiv.org/abs/1903.10456} {\path{arXiv:1903.10456}},
  \href {https://doi.org/10.3847/1538-4365/ab18ff}
  {\path{doi:10.3847/1538-4365/ab18ff}}.

\bibitem{Boscarino.2015}
S.~Boscarino, R.~B\"urger, P.~Mulet, G.~Russo, L.~M. Villada, {Linearly
  Implicit IMEX Runge--Kutta Methods for a Class of Degenerate
  Convection-Diffusion Problems}, SIAM Journal on Scientific Computing 37~(2)
  (2015) B305--B331.
\newblock \href {https://doi.org/10.1137/140967544}
  {\path{doi:10.1137/140967544}}.

\bibitem{Boscarino.2016}
S.~Boscarino, F.~Filbet, G.~Russo, {High Order Semi-implicit Schemes for Time
  Dependent Partial Differential Equations}, Journal of Scientific Computing
  68~(3) (2016) 975--1001.
\newblock \href {https://doi.org/10.1007/s10915-016-0168-y}
  {\path{doi:10.1007/s10915-016-0168-y}}.

\bibitem{Boscarino.20168e8}
S.~Boscarino, R.~Bürger, P.~Mulet, G.~Russo, L.~M. Villada, {On linearly
  implicit IMEX Runge-Kutta methods for degenerate convection-diffusion
  problems modeling polydisperse sedimentation}, Bulletin of the Brazilian
  Mathematical Society, New Series 47~(1) (2016) 171--185.
\newblock \href {https://doi.org/10.1007/s00574-016-0130-5}
  {\path{doi:10.1007/s00574-016-0130-5}}.

\bibitem{boscarino2018all}
S.~Boscarino, G.~Russo, L.~Scandurra, All mach number second order
  semi-implicit scheme for the euler equations of gas dynamics, Journal of
  Scientific Computing 77 (2018) 850--884.

\bibitem{avgerinos2019linearly}
S.~Avgerinos, F.~Bernard, A.~Iollo, G.~Russo, Linearly implicit all mach number
  shock capturing schemes for the euler equations, Journal of Computational
  Physics 393 (2019) 278--312.

\bibitem{Gonzalez-Pinto.2022}
S.~González-Pinto, D.~Hernández-Abreu, M.~S. Pérez-Rodríguez, A.~Sarshar,
  S.~Roberts, A.~Sandu, {A unified formulation of splitting-based implicit time
  integration schemes}, Journal of Computational Physics 448 (2022) 110766.
\newblock \href {http://arxiv.org/abs/2103.00757} {\path{arXiv:2103.00757}},
  \href {https://doi.org/10.1016/j.jcp.2021.110766}
  {\path{doi:10.1016/j.jcp.2021.110766}}.

\bibitem{Sandu.2015}
A.~Sandu, M.~Günther, {A Generalized-Structure Approach to Additive
  Runge--Kutta Methods}, SIAM Journal on Numerical Analysis 53~(1) (2015)
  17--42.
\newblock \href {https://doi.org/10.1137/130943224}
  {\path{doi:10.1137/130943224}}.

\bibitem{cooper1980additive}
G.~J. Cooper, A.~Sayfy, Additive methods for the numerical solution of ordinary
  differential equations, Mathematics of Computation 35~(152) (1980)
  1159--1172.

\bibitem{Cooper.1983sej}
G.~J. Cooper, A.~Sayfy, {Additive Runge-Kutta methods for stiff ordinary
  differential equations}, Mathematics of Computation 40~(161) (1983) 207--218.
\newblock \href {https://doi.org/10.1090/s0025-5718-1983-0679441-1}
  {\path{doi:10.1090/s0025-5718-1983-0679441-1}}.

\bibitem{Ascher.1997}
U.~M. Ascher, S.~J. Ruuth, R.~J. Spiteri, {Implicit-explicit Runge-Kutta
  methods for time-dependent partial differential equations}, Applied Numerical
  Mathematics 25~(2-3) (1997) 151--167.
\newblock \href {https://doi.org/10.1016/s0168-9274(97)00056-1}
  {\path{doi:10.1016/s0168-9274(97)00056-1}}.

\bibitem{Lowrie.2008}
R.~B. Lowrie, J.~D. Edwards, {Radiative shock solutions with grey
  nonequilibrium diffusion.}, Shock Waves 18 (2008) 129--143.

\bibitem{Pareschi.2005}
L.~Pareschi, G.~Russo, {Implicit–Explicit Runge–Kutta Schemes and
  Applications to Hyperbolic Systems with Relaxation}, Journal of Scientific
  Computing 25~(1) (2005) 129--155.
\newblock \href {https://doi.org/10.1007/s10915-004-4636-4}
  {\path{doi:10.1007/s10915-004-4636-4}}.

\bibitem{Giraldo.2013}
F.~X. Giraldo, J.~F. Kelly, E.~M. Constantinescu, {Implicit-Explicit
  Formulations of a Three-Dimensional Nonhydrostatic Unified Model of the
  Atmosphere (NUMA)}, SIAM Journal on Scientific Computing 35~(5) (2013)
  B1162--B1194.
\newblock \href {https://doi.org/10.1137/120876034}
  {\path{doi:10.1137/120876034}}.

\bibitem{Conde.2017}
S.~Conde, S.~Gottlieb, Z.~J. Grant, J.~N. Shadid, {Implicit and
  Implicit–Explicit Strong Stability Preserving Runge–Kutta Methods with
  High Linear Order}, Journal of Scientific Computing 73~(2-3) (2017) 667--690.
\newblock \href {https://doi.org/10.1007/s10915-017-0560-2}
  {\path{doi:10.1007/s10915-017-0560-2}}.

\bibitem{Sebastiano.2023}
S.~Boscarino, {High-Order Semi-implicit Schemes for Evolutionary Partial
  Differential Equations with Higher Order Derivatives}, Journal of Scientific
  Computing 96~(1) (2023) 11.
\newblock \href {https://doi.org/10.1007/s10915-023-02235-0}
  {\path{doi:10.1007/s10915-023-02235-0}}.

\bibitem{imex-trt}
B.~S. Southworth, S.~S. Olivier, H.~Park, T.~Buvoli, One-sweep moment-based
  semi-implicit-explicit integration for gray thermal radiation transport,
  arXiv preprint arXiv:2401.04285 (2024).

\bibitem{Soderlind.2006}
G.~S\"oderlind, L.~Wang, {Adaptive time-stepping and computational stability},
  Journal of Computational and Applied Mathematics 185~(2) (2006) 225--243.
\newblock \href {https://doi.org/10.1016/j.cam.2005.03.008}
  {\path{doi:10.1016/j.cam.2005.03.008}}.

\bibitem{Ranocha.2022}
H.~Ranocha, L.~Dalcin, M.~Parsani, D.~I. Ketcheson, {Optimized Runge-Kutta
  Methods with Automatic Step Size Control for Compressible Computational Fluid
  Dynamics}, Communications on Applied Mathematics and Computation 4~(4) (2022)
  1191--1228.
\newblock \href {http://arxiv.org/abs/2104.06836} {\path{arXiv:2104.06836}},
  \href {https://doi.org/10.1007/s42967-021-00159-w}
  {\path{doi:10.1007/s42967-021-00159-w}}.

\bibitem{nprk1}
T.~Buvoli, B.~S. Southworth, A new class of runge-kutta methods for nonlinearly
  partitioned systems, arXiv preprint arXiv:2401.04859 (2024).

\bibitem{nprk2}
B.~K. Tran, B.~S. Southworth, T.~Buvoli, Order conditions for nonlinearly
  partitioned runge-kutta methods, arXiv preprint arXiv:2401.12427 (2024).

\end{thebibliography}
\end{document}